\documentclass{article} 

\usepackage{scalefnt}
\usepackage{amsmath,amssymb}
\usepackage{bbm}
\usepackage{wrapfig}
\usepackage{verbatim}
\usepackage{epsfig}
\usepackage{textcomp}
\usepackage{vmargin}

\usepackage{graphicx}          
 
\usepackage{scalerel}
\usepackage{stackengine,wasysym}

\newcommand\reallywidetilde[1]{\ThisStyle{%
  \setbox0=\hbox{$\SavedStyle#1$}%
  \stackengine{-.1\LMpt}{$\SavedStyle#1$}{%
    \stretchto{\scaleto{\SavedStyle\mkern.2mu\AC}{.5150\wd0}}{.6\ht0}%
  }{O}{c}{F}{T}{S}%
}}

\allowdisplaybreaks

\def\d{\displaystyle}
\def\R{\mathbb{R}}

\def\E{\mathbb{E}}
\def\P{\mathbb{P}}

\def\K{\mathbb{K}}
\def\X{\mathbb{X}}
\def\Y{\mathbb{Y}}

\newcommand{\fF}{f\!F}

\newcommand{\fFbar}{{\overline{f\!F}}}
\newcommand{\fFtilde}{{\reallywidetilde{f\!F}}}
\newcommand{\fFbarp}{{\overline{f\!F}^+}}

\newcommand{\thetabar}{{\bar{\theta}}}
\newcommand{\Thetabar}{{\bar{\Theta}}}

\newcommand{\Psibar}{{\bar{\Psi}}}
\newcommand{\Rbar}{{\bar{R}}}

\newcommand{\ybar}{{\bar{y}}}
\newcommand{\Ybar}{{\bar{Y}}}

\newcommand{\vbar}{{\bar{v}}}

\newcommand{\thetahat}{{\hat{\theta}}}
\newcommand{\Thetahat}{{\hat{\Theta}}}

\newcommand{\Rhat}{{\hat{R}}}

\newcommand{\yhat}{{\hat{y}}}
\newcommand{\Yhat}{{\hat{Y}}}

\newcommand{\vhat}{{\hat{v}}}

\newcommand{\Pcal}{\mathcal P}

\def\Mtn1{M_{t_{n+1}}}
\def\Mtn{M_{t_{n}}}
\def\Xtn1{X_{t_{n+1}}}
\def\Xtn{X_{t_{n}}}
\def\Ytn1{Y_{t_{n+1}}}
\def\Ytn{Y_{t_{n}}}
\def\phiilip{[\Phi_i]}
\def\phi0lip{[\Phi_0]}
\def\philip{[\Phi]}

\newcommand{\barklip}[1]{[#1]_{{k}}}
\newcommand{\barkplip}[1]{[#1]_{{k+1}}}

\newcommand{\barnlip}[1]{[#1]_{{n}}}
\newcommand{\barnplip}[1]{[#1]_{{n+1}}}

\newcommand{\nsup}[1]{\|#1\|}
\newcommand{\barksup}[1]{\|#1\|_{{k}}}
\newcommand{\barkpsup}[1]{\|#1\|_{{k+1}}}

\newcommand{\barnpsup}[1]{\|#1\|_{{n+1}}}
\newcommand{\dk}[1]{\|#1\|_{{k,1}}}
\newcommand{\dkp}[1]{\|#1\|_{{k+1,1}}}

\newcommand{\dn}[1]{\|#1\|_{{n,1}}}

\newcommand{\1}{\mathbbm{1}}

\usepackage[]{xcolor}
\definecolor{darkred}{rgb}{0.65,0.15,0.25}
\definecolor{darkgreen}{rgb}{0,0.4,0}
\definecolor{myblue}{rgb}{0.2,0.2,0.7}
\definecolor{aliceblue}{rgb}{0.94, 0.97, 1.0}
\definecolor{airforceblue}{rgb}{0.36, 0.54, 0.66}
\definecolor{antiquebrass}{rgb}{0.8, 0.58, 0.46}
\definecolor{cambridgeblue}{rgb}{0.64, 0.76, 0.68}
\definecolor{asparagus}{rgb}{0.53, 0.66, 0.42}
\definecolor{chr1}{HTML}{A6CEE3}
\definecolor{chr2}{HTML}{1F78B4}
\definecolor{chr3}{HTML}{B2DF8A}
\definecolor{chr4}{HTML}{33A02C}
\definecolor{chr5}{HTML}{FB9A99}
\definecolor{chr6}{HTML}{E31A1C}
\definecolor{chr7}{HTML}{FDBF6F}
\definecolor{chr8}{HTML}{FF7F00}
\definecolor{chr9}{HTML}{CAB2D6}
\definecolor{chr10}{HTML}{6A3D9A}
\definecolor{chr11}{HTML}{FFFF99}
\definecolor{chr12}{HTML}{B15928}
\definecolor{chr13}{HTML}{8DD3C7}
\definecolor{chr14}{HTML}{FFFFB3}
\definecolor{chr15}{HTML}{BEBADA}
\definecolor{chr16}{HTML}{FB8072}
\definecolor{chr17}{HTML}{80B1D3}
\definecolor{chr18}{HTML}{FDB462}
\definecolor{chr19}{HTML}{B3DE69}
\definecolor{chr20}{HTML}{FCCDE5}
\definecolor{chr21}{HTML}{D9D9D9}
\definecolor{chr22}{HTML}{BC80BD}
\definecolor{chrX}{HTML}{CCEBC5}
\definecolor{chrY}{HTML}{FFED6F}

\usepackage{tikz}
\usetikzlibrary{decorations.markings}
\usetikzlibrary{decorations.pathreplacing}
\usetikzlibrary{spy}
\usetikzlibrary{positioning}
\usetikzlibrary{backgrounds}
\usetikzlibrary{hobby,decorations}
\usetikzlibrary{shapes,snakes}
\usetikzlibrary{shapes.geometric}
\usetikzlibrary{shadows}
\usetikzlibrary{intersections}
\tikzset{   invisible/.style={opacity=0,text opacity=0},     visible on/.style={alt={#1{}{invisible}}},     alt/.code args={<#1>#2#3}{%
     \alt<#1>{\pgfkeysalso{#2}}{\pgfkeysalso{#3}} },}
\tikzset{cross/.style={cross out, draw=black, minimum size=2*(#1-\pgflinewidth), inner sep=0pt, outer sep=0pt},cross/.default={1pt}}

\def\newblock{\hskip .11em plus .33em minus .07em}

\newtheorem{theor}{Theorem}[section]
\newtheorem{lemma}[theor]{Lemma}
\newtheorem{propo}[theor]{Proposition}
\newtheorem{defnt}[theor]{Definition}

\newtheorem{exam}{Example}

\newcounter{sousexam}[exam]
\newenvironment{ssexam}
{\par \refstepcounter{sousexam} Model \arabic{exam}.\alph{sousexam}-- \rm}
{\par\ifdim\lastskip<3pt \removelastskip \penalty-200 \vskip3pt \fi}
\newcommand{\bs}[1]{\boldsymbol{#1}}

\DeclareMathAlphabet\mathbfcal{OMS}{cmsy}{b}{n}

\begin{document}

\title{Change-point detection for Piecewise Deterministic Markov Processes}                                              

\author{Alice Cleynen \and Beno\^\i te de Saporta\thanks{                             
The work was partially supported by R\'egion Languedoc-Roussillon and FEDER under grant \emph{Chercheur(se)s d'Avenir}, project PROMMECE.}}       
\date{\small Institut Montpelli\'erain Alexander Grothendieck, CNRS, Univ. Montpellier, France}

\maketitle                       

\begin{abstract}                          
We consider a change-point detection problem for a simple class of Piecewise Deterministic Markov Processes (PDMPs). A continuous-time PDMP is observed in discrete time and through noise, and the aim is to propose a numerical method to accurately detect both the date of the change of dynamics and the new regime after the change. To do so, we state the problem as an optimal stopping problem for a partially observed discrete-time Markov decision process taking values in a continuous state space and provide a discretization of the state space based on quantization to approximate the value function and build a tractable stopping policy. We provide error bounds for the approximation of the value function and numerical simulations to assess the performance of our candidate policy.
\end{abstract}
%
\section{Introduction}
\label{Introduction}
%
Piecewise Deterministic Markov processes (PDMPs) are a general class of non-diffusion processes introduced by M. Davis in the 80's \cite{Davis84} covering a wide range of applications from workshop optimization, queuing theory \cite{Davis93}, internet networks \cite{BCGMZ13}, reliability \cite{dSDZ16}, insurance and finance \cite{bauerle11} or biology \cite{DHKL15,RT15,RTW12} for instance. 
PDMPs are continuous time hybrid processes with a discrete component called mode or regime and a Euclidean component. The process follows a deterministic trajectory punctuated by random jumps. In the special case where the Euclidean component is continuous the jumps correspond to a change of regime. For many applications, the regime is not observed and the Euclidean variable is measured in discrete-time, through noise. It may be e.g. a degradation or failure of some component of a system, see \cite{BBGPPS14} where the Euclidean component is some cool down time that increases with the degradation of the system, or the cancer cell load of remission patients monitored through proxy tumor markers at regular follow-up blood tests to detect relapse \cite{AM06}.
The aim of this paper is to propose a fully computable discretization of the value function of the optimal stopping problem corresponding to the change-point detection, and derive error bounds for this approximation.  We also use the approximation to build a computable candidate strategy that should be close to optimality. We assess its performance on numerical examples.

The general problem of change-point detection can be seen as an impulse control problem if there are multiple changes in regime. This is a very difficult problem.
Although the optimal control of PDMPs has attracted a lot of attention since the 80s, see e.g. \cite{Elliott12,CD13,Davis93,Dempster95,Gatarek92,Lenhart89}, very few works consider such models under partial observations.
In \cite{BdSD13}, the authors consider an optimal stopping problem for PDMPs where the jump times are perfectly observed and the post-jump locations are observed through noise. They derive the dynamic programming equations of the problem, as well as a numerical approximation of the value function and a computable $\epsilon$-optimal stopping time.
In \cite{BL17}, the authors consider a general continuous control problem where both the jump times and post-jump locations are observed through noise. They  reduce the problem to a discrete-time Markov Decision Process (MDP) and prove the existence of optimal policies, but provide no numerical approximation of the value function or optimal strategies.

In the present paper, we make a first step towards solving the difficult problem of change-point detection of PDMPs when the jumps are not observed at all. We address the simple case where there is only one change of regime to detect. The problem can thus be formulated as an optimal stopping problem for PDMPs under partial observations. However, unlike \cite{BL17,BdSD13} we do not suppose that the observations are made at or around the jump times. Instead, we suppose that the observations times are deterministic and on a regular grid of step size $\delta$. This enables us to formulate the problem as an optimal stopping problem for a discrete-time partially observed MDP. The equivalent fully-observed MDP for the filter process is still in discrete-time but on an infinite state space. We then propose a two-step discretization of this MDP, following an idea from \cite{dSDN16}. The first step is a time-dependent discretization of the state space of the original PDMP. The second step is a joint discretization of the approximate filter thus obtained together with an approximation of the observation process. Note that unlike \cite{dSDN16} or \cite{bauerle11}, we do not make the assumption that the MDP kernel has a density with respect to some fixed probability measure. 

The main discretization tool we use is optimal quantization. The quantization of a random variable $X$ consists in finding a finite grid such that the projection $\widehat{X}$ of $X$ on this grid minimizes some $L_{p}$ norm of the difference $X-\widehat{X}$. 
There exists an extensive literature on quantization methods for random variables and processes. 
The interested reader may for instance consult \cite{gray98,pages04} and the references therein. 
Quantization methods have been developed recently in numerical probability or optimal stochastic control with applications in finance, see e.g. \cite{bally03,bally05,pages98,pages04}. 

The paper is organized as follows. In Section \ref{Framework}, we introduce our continuous-time PDMP model as well as the observation model. We define the change-point detection problem as an optimal stopping problem under partial observations and give the equivalent fully observed dynamic programming equations for the filter process. In Section \ref{sec:value}, we propose a two-step discretization approach by quantization to numerically solve the optimization problem and build a tractable strategy. Proofs of our main statements are postponed to Section \ref{sec:proofs}. In Section \ref{Simulations}, we investigate the performance of our candidate strategy and compare our approach to moving average and Kalman filtering when possible. A conclusion is given in Section \ref{sec:CCL} and the proofs are gathered in the appendix.
%
\section{Model and problem setting} 
\label{Framework}
%
In this section, we present the special class of PDMPs we focus on, define the observation process and state the change-point detection problem as an optimal stopping problem under partial observation. We then derive the filter recursive equation and state the equivalent fully observed optimal stopping problem as well as the corresponding dynamic programming equations.
%
\subsection{Continuous-time PDMP model}
%
We consider the problem of detecting a change-point in the dynamic of a special class of PDMPs which is observed with noise on discrete observation times. 
The process  $\bs{X}_t=(m_t,x_t,t)$ is defined on a state space $\bs{E}=\mathcal{M}\times\K\times\mathbb{R}_+$, where $\mathcal{M}=\{0,1,\dots,d\}$ is the finite set of modes and $\K$ is a compact subset of $\R$. We will further denote $\X=\mathcal{M}\times \K$. For each mode $m$, the local characteristics of the PDMP are\\
$\bullet$ a flow $\bs{\Psi}_m:\K\times\mathbb{R}_+^{2}\rightarrow\K\times\mathbb{R}_+$ of the special form
$\bs{\Psi}_m((x,u),t)=(\Phi_m(x,t),u+t)$;\\
$\bullet$ a jump intensity $\lambda_m:\K\times\mathbb{R}_+\rightarrow\mathbb{R}_+$ such that $\lambda_0(x,u)=\lambda(u)$ and $\lambda_i(x,u)=0$, for all positive $i$ in $\mathcal{M}$;\\
$\bullet$ a (sub)Markovian jump kernel $\bs{Q}_m:(\mathcal{B}(\bs{E}),\K\times\mathbb{R}_+)\rightarrow\bs{E}$ such that
$\bs{Q}_m(\{i\}\times A\times\{0\}|x,u) =\pi_i\1_{\{0\}}(m)\1_{A}(x)$ with $\pi_i>0$ for all positive $i\in \mathcal{M}$ and $\sum_{i=1}^d \pi_i = 1$.

In other words, the PDMP has a single jump at some random time $T$ and evolves deterministically after the jump. The distribution of the jump time does not depend on the position $x_t$ but only on the running time:
\begin{equation*}
\mathbb{P}_{(0,x,0)}(T>t)=e^{-\int_0^t \lambda(s)ds}.
\end{equation*}
After the jump, the location $x$ remains unchanged, the time since the last jump $u$ is set to $0$ and a new mode is selected according to the distribution $\pi$. The third component of $\bs{X}_t$, namely the running time since the last jump, only intervenes in the jump time distribution. It  is necessary to obtain a strong Markov process. 
The assumption that the flow $\Phi_m$ does not depend on the running time $u$ is made only to keep notation simple and is not actually required, see Example~\ref{eq:exsin2}.
As we will see in the sequel, solving the change-point detection problem is not straightforward, even for such simple dynamics.

We suppose that the observation times $(t_n)_{n\in\mathbb{N}}$ are deterministic and on a regular grid of step size $\delta$ until a finite horizon $N\delta$, and that a noisy observation of $x_{t_n}$ is available at each time $t_n$:
\begin{equation}
\label{Obs}
{Y}_n=F(\bs{X}_{t_n})+\varepsilon_n=F(x_{t_n})+\varepsilon_n, 
\end{equation} 
where $F$ is a deterministic link function, $(\varepsilon_n)$ are \textit{iid} real-valued random variables with density $f$ with respect to the Lebesgue measure on $\R$ and independent from the process $(\bs{X}_t)$. We further assume that $Y$ takes its values in $\Y$, subset of $\R$. 
%
\subsection{Examples}
\label{sec:examples}
%
The following toy examples will be extensively investigated numerically in Section~\ref{Simulations}. In all examples, the jump intensity is of the form $\lambda(t)=t$ so that the probability to jump gets higher as time goes by. The Markov kernel $Q$ is the uniform distribution on the possible post-jump modes. 
The distribution of the noise is a centered Gaussian with variance $\sigma^2$ truncated at $[-s,s]$ for some $s\in \R$.
We investigate several forms for the flow. The link function between the process and the observations will be either $F(x)=x$ or $F(x)=x^{-1}$.
%
\begin{exam} \label{Ex:expo}
In order to compare our method with other state of the art approaches, we study exponential or linear trajectories. The process starts at $\bs{X}_0= (0,1,0)$. 
\begin{ssexam} The flows $\Phi_m$ are defined as
\begin{equation*} \label{eq:exexpo}
  \left\{
  \begin{array}{rl}
\Phi_0(x,t)&=x, \\
\Phi_i(x,t)&=e^{v_it}x, \quad v_i\neq 0; 1\leq i\leq d.
 \end{array}
  \right.
\end{equation*}
\end{ssexam}
\begin{ssexam} The flows $\Phi_m$ are defined as
\begin{equation*} \label{eq:exexpob}
  \left\{
  \begin{array}{rl}
\Phi_0(x,t)&=x, \\
\Phi_1(x,t)&=e^{v_1t}x, \quad v_1\neq 0,\\
\Phi_2(x,t)&=x+v_2t, \quad v_2\neq 0.
 \end{array}
  \right.
\end{equation*}
\end{ssexam}
\end{exam}
%
\begin{exam} \label{Ex:sin1}
We study the more challenging example of non one-to-one flows, trying to detect a sudden change in either the frequency or the slope of a sinusoidal trajectory. The process is initiated  at $\bs{X}_0= (0,0,0)$.
\begin{ssexam} The flows $\Phi_m$ are defined as
\begin{equation*}\label{eq:exsin1}
  \left\{
  \begin{array}{rl}
 \Phi_0(x,t)&=\sin(\arcsin(x)+v_0\pi t), \\
 \Phi_i(x,t)&=\sin(\arcsin(x)+v_i\pi t), 
 \end{array}
  \right.
\end{equation*}
with $v_i\neq v_0$ for $1\leq i\leq d$.
\end{ssexam}
\begin{ssexam} The flows $\bs\Psi_m$ are defined as
\begin{equation*} \label{eq:exsin2}
  \left\{
  \begin{array}{rl}
 \bs\Psi_0(x,u,t)&=(\sin(\arcsin(x)+v_0\pi t),u+t),\\
 \bs\Psi_i(x,u,t)&=(\sin(\arcsin(x- v_iu)+v_0\pi t) + v_i (u+t),\\&\quad u+t).
 \end{array}
  \right.
\end{equation*}
with $v_i\neq 0$ for $1\leq i\leq d$.
\end{ssexam}
\end{exam}
%
\subsection{Optimal stopping problem under partial observations}
\label{sec:part obs}
%
We are interested in detecting the jump-time $T$ and the mode after the jump based on the observations $Y_n$. We choose to formulate this problem as an optimal stopping problem for a discrete-time Markov chain. However in our framework, it is important to note that the underlying process is time-continuous, and in particular that the jump-time $T$ may occur between observation dates. In this paper, we will only allow detections at the observation times. Allowing detection between observation times is still an open problem.

In the sequel, we will simply denote $X_n=(m_{t_n},x_{t_n})$. As our PDMP has only one jump, one can explicitly write the kernels $P_n$ of the time-inhomogeneous discrete-time Markov chain $(X_n)$.  For any Borelian subsets $A\subset\mathcal{M}$, $B \subset\K$, any $(m,x) \in \X$ and $n\geq 0$, one has
\begin{align}
{P}_n(A\times B|m,x) \label{def:Pn}
  &=\mathbb{P}({X}_{{n+1}}\in A\times B| {X}_{n}=(m,x))\\
  &=\1_{\{0\}}(m)\1_{A}(0)\1_{B}(\Phi_0(x,\delta))e^{-\int_0^\delta \lambda(n\delta+s)ds}\nonumber\\
  &+\sum_{i=1}^d\pi_i\1_{\{0\}}(m)\1_{A}(i)\int_0^\delta \lambda(n\delta+s)e^{-\int_0^s \lambda(n\delta+z)dz}\1_{B}(\Phi_i(\Phi_0(x,s),\delta-s))ds\nonumber\\
&+\sum_{i=1}^d\1_{\{i\}}(m)\1_{A}(i)\1_{B}(\Phi_i(x,\delta)).\nonumber
\end{align}
For any $y'\in\Y$, let $f\!F_{y'}$ be the function from $\X$ onto $\R$ defined by 
\begin{equation*}
f\!F_{y'}:(m,x)\mapsto f(y'-F(x))
\end{equation*}
Thus the kernels $R_n$ of the Markov chain $(X_n,Y_n)$ are, for any Borelian subsets $A\subset\X$, $C\subset\Y$, any $(m,x,y) \in \X\times\Y$ and $n\geq 0$
\begin{align}
{R}_n(A\times C|m,x,y)\label{def:Rn}
  &=\mathbb{P}(({X}_{{n+1}},{Y}_{n+1})\in A\times C| ({X}_{n},{Y}_n)=(m,x,y))\\
  &=\int_C P_n(f\!F_{y'}\1_A)(m,x)dy',\nonumber
\end{align}
Note that the $R_n$ kernels do not depend on $y$. 

We can now state our change-point detection problem. For $0\leq n\leq N$, set $\mathcal{F}_n=\sigma(X_{k},Y_k,0\leq k\leq n)$ the $\sigma$-field generated by the Markov chain $(X_n,Y_n)$ up to time $n$, and
$\mathcal{F}_n^\Y=\sigma(Y_k,0\leq k\leq n)$ the $\sigma$-field generated by the observations up to time $n$.
Let $\mathcal{T}^\mathbb{Y}$ be the set of $\mathcal{F}^\Y$-stopping times. We do not restrict ourselves to the stopping times bounded by $N$ because it may be optimal not to stop at all during until the horizon $N$ is reached.
A decision taken at the stopping-time $\tau\in\mathcal{T}^\mathbb{Y}$ is a $\mathcal{F}_{\tau}^\Y$-measurable random variable $A$ taking values in $\mathcal{M}_+=\{1,2,\ldots,d\}$ if $\tau \leq N$, equal to $0$ if $\tau>N$. Decision $A=i$ corresponds to deciding $X_{\tau}$ is in mode $i$. Let $\mathcal{A}_\tau^\mathbb{Y}$ be the set of admissible decisions at stopping time $\tau$.

Until stopping-time $\tau$, the cost-per-stage function is denoted by $c$ and the terminal cost (at stopping time $\tau$) when taking decision $a$ is $C$, where
\begin{eqnarray*}
c(i,x,y)&=&\beta_i\delta, \\
C(m,x,y,0)&=&c(m,x,y), \\
C(m,x,y,a>0)&=&\alpha\1_{(m=0)} +  \gamma_{ma}  \1_{(m\neq a; m>0)}, 
\end{eqnarray*}
with $\beta_0=0$ and for positive $i$, $\beta_i=\beta>0$. Thus, $\beta$ represents the penalty for late detection, $\alpha$ the false alarm penalty and $\gamma_{ma}$ the cost of selecting mode $a$ instead of mode $m$. The cost of an admissible strategy $(\tau,A_\tau)\in\mathcal{T}^\mathbb{Y}\times \mathcal{A}_\tau^\mathbb{Y}$ for starting point $\xi\in\X\times\Y$ is
\begin{align*}
J(\tau,A_\tau,\xi) 
  &={\mathbb{E}}_{\xi}\Bigg[\!\sum_{n=0}^{(\tau-1)\wedge N}\! c(X_{n},Y_n)\!+\!\1_{(\tau\leq N)}C(X_{\tau\wedge N},Y_{\tau\wedge N},A_\tau)\Bigg],
\end{align*}
and the value function of the problem is
\begin{equation*}
V(\xi)=\inf_{(\tau,A_\tau) \in \mathcal{T}^\mathbb{Y}\times \mathcal{A}_\tau^\Y} J(\tau,A_\tau,\xi').
\end{equation*} 
%
The optimal (possibly not achievable) cost is $0$ when the jump is detected at the first observation after its occurrence and the right post-jump mode is selected. 
The aim of this paper is to derive a numerically tractable approximation of the value function $V$ as well as propose a computable strategy close to optimality.
%
\subsection{Fully observed optimal stopping problem}
The classical approach to deal with partial observations is to introduce the filter process and the corresponding completely observed optimal stopping problem for filtered trajectories.
For any starting point $\xi=(0,x,y)\in\X\times\Y$, set $\Theta_0=\theta_0=\delta_{(0,x)}$ and for $1\leq n\leq N$, and any Borelian subset $A$ of $\X$ set
\begin{eqnarray*}
\Theta_n(A)&=&\mathbb{P}_{\xi}(X_{n}\in A| \mathcal{F}_n^\mathbb{Y}),
\end{eqnarray*}
 the filter for the unobserved part of the  process. The filter is recursively obtained as follows.
%
\begin{propo}\label{def:psi} For any $n\geq 0$, conditionally on $(\Theta_{n}=\theta,Y_{n+1}=y')$, one has $\Theta_{n+1}=\Psi_{n}(\theta,y')$ with 
\vspace{-.2cm}
\begin{align}
\Psi_{n}(\theta,y')(A)
&= \frac{\int_{\mathbb{X}}P_n(f\!F_{y'}\1_A)(m,x)d\theta(m,x)}{\int_{\mathbb{X}}P_n(f\!F_{y'})(m,x)d\theta(m,x)},
\end{align} \label{Psi}
for any Borelian subset $A$ of $\X$. 
\end{propo}\label{correction}
%
The proof of this proposition is quite classical and therefore omitted. It relies on the standard prediction-correction approach. Similar computations can be found e.g. in \cite{bauerle11} in the framework of MDPs, with the notable difference that in our context we do not assume that the kernels have a density with respect to any fixed measure, or in \cite{BdSD13} for a different class of PDMPs.

Set $\mathcal{P}(\X)$ the set of probability measures on $\X$.
Thus, $(\Theta_n,Y_n)$ is a Markov chain on $\mathcal{P}(\X)\times\Y$, with transition kernels defined, for any Borelian subsets $P\subset\mathcal{P}(\X)$, $C\subset\Y$, and any $(\theta,y)\in \Pcal(\X)\times\Y$, by
\begin{align*}
{R'}_n(P\times C|\theta,y)
   &=\mathbb{P}(({\Theta}_{{n+1}},Y_{n+1})\in P\times C| ({\Theta}_{n},Y_n)=(\theta,y))\\
&=\int_{\X\times C}\1_{P}\left(\Psi_n(\theta,y') \right)P_n(f\!F_{y'})(m,x)d\theta(m,x)dy'.
\end{align*}
%
Again, this kernel does not depend on $y$.
The partially observed optimal stopping problem defined in Section~\ref{sec:part obs} is equivalent to a fully observed optimal stopping problem using the filtered trajectories introduced above. The fully observed state space is thus $\mathcal{P}(\X)\times\Y$, the initial point is $\xi'=(\delta_{(0,x)},y)$ for some $(x,y)\in\K\times\Y$. 
In this framework, the cost of an admissible strategy $(\tau,A_\tau)\in\mathcal{T}^\mathbb{Y}\times \mathcal{A}_\tau^\mathbb{Y}$ for starting point $\xi\in\X\times\Y$ is
\begin{align*}
J'(\tau,A_\tau,\xi') 
  &={\mathbb{E}}_{\xi'}\!\Bigg[\sum_{n=0}^{(\tau-1)\wedge N}\! c'(\Theta_{n},Y_n)\!+\!\1_{(\tau\leq N)}C'(\Theta_{\tau\wedge N},Y_{\tau\wedge N},A_\tau)\Bigg],
\end{align*}
where, for $g$ from $\X\times\Y$ onto $\R$  $g'$ is the function from $\mathcal{P}(\X)\times\Y$ onto $\R$ such that $g'(\theta,y)=\int_\X g(m,x,y)d\theta(m,x)$, here for  $g=c$ or $g=C$.
The value function of the problem is
\begin{equation}
V'(\xi')=\inf_{(\tau,A_\tau) \in \mathcal{T}^\mathbb{Y}\times \mathcal{A}_\tau^\Y} J'(\tau,A_\tau,\xi'). \label{eq:defV}
\end{equation} 
The value function is then solution to the dynamic programming equations. 
\begin{theor}\label{th:DP'} Set $v_N'(\theta,y)= \min_{a\in\mathcal{M}} \ C'(\theta, y,a)$ and for $0\leq n\leq N-1$
    \begin{eqnarray*}
      v_n'(\theta,y)&=&\min\left\{\min_{a\in\mathcal{M}_+} C'(\theta,y,a); c'(\theta,y)+\R'_nv'_{n+1}(\theta,y) \right\}.
    \end{eqnarray*}
Let $\xi_0'=(\delta_{(0,x)},y) \in \mathcal{P}(\X)\times \Y$. Then we have $$v'_0(\xi_0') = V'(\xi_0')=V(0,x,y).$$ 
\end{theor}\label{dyn-prog}
%
Note that none of the functions above actually depends on $y$.
Again, the proof of this statement relies on standard arguments and is omitted. The proper framework for the proof is that of Partially Observed MDPs (POMDPs). One first defines the equivalent POMDP to the optimal stopping problem under partial observation and then proves the equivalence with the fully observed MDP corresponding to the fully observed optimal stopping problem. The dynamic programming is then straightforward. Similar derivations can be found for instance in \cite{bauerle11,dSDN16}.
%
\section{Numerical approximation of the value functions}
\label{sec:value}
The aim of this paper is to propose a numerically tractable approximation of the optimal value function $V'$ defined in eq.~(\ref{eq:defV}) and a corresponding candidate optimal strategy. 
The main difficulties are first that the filter $\Theta_n$ is measure-valued and thus infinite-dimensional and second that this filter cannot be simulated as the Bayes operators $\Psi_n$ involve continuous integration.

\begin{figure}[h!]
  \begin{center}
    \framebox{
      \begin{tikzpicture}[,scale=1, every node/.style={scale=0.7}]
      \node(cont) at (-1, 2.5) {$\bs{X}_t=(m_t,x_t,t)$};
      \node(cont-EP) at (-1, 2) {$\X\times\R^+$, $\bs{P}$};
      \node(in) at (-1,0.5) {${X}_{n}=(m_{t_n},x_{t_n})$};
      \node(in-EP) at (-1, 0) {$\X$, ${P}_n$};
      \draw[->] (cont-EP)--(in);
      \node(va) at (-1,-1.5) {$({X}_{n},{Y}_n)$};
      \node(va-EP) at (-1, -2) {$\X\times\Y$, ${R}_n$};
      \draw[->] (in-EP) --(va) node[left=0.15cm,pos=0.5,color=blue!50] {\small observations} node[right=0.15cm,pos=0.5,color=chr16] {\tiny ${Y}_n=F({X}_{n})+\varepsilon_n$};
      \node(fi) at (-1,-3.5) {$(\Theta_n,Y_n)$};
      \node(fi-EP) at (-1, -4) {$\mathcal{P}(\X)\times\Y$, ${R}'_n$};
      \node(dyn) at (-1,-5.5) {$v'_n(\Theta_n,Y_n)$};
      \draw[->] (va-EP)--(fi) node[left=0.15cm,pos=0.5,color=blue!50] {\small filtering} node[right=0.15cm,pos=0.5,color=chr16] {\small $\Psi$};
      \draw[->] (fi-EP)--(dyn) node[left=0.15cm,pos=0.35,color=blue!50] {\small dynamic} node[left=0.15cm,pos=0.65,color=blue!50] {\small programming}; 
      \node(q1) at (3,0.5) {$(\bar{m}_{t_n},\bar{x}_{t_n})=\bar{X}_{n}$};
      \node(q1-EP) at (3,0) {$\Omega_n$, $\bar{{P}}_n$};%
      \node(vaq1) at (3,-1.5) {$(\bar{X}_{n},\Ybar_n)$};
      \node(vaq1-EP) at (3, -2) {$\Omega_n\times\Y$, $\bar{R}_n$};
      \node(fiq1) at (3,-3.5) {$(\bar{\Theta}_n,\Ybar_n)$};
       \node(fiq1-EP) at (3, -4) {$\mathcal{P}(\Omega_n)\times\Y$, $\bar{R}'_n$};
      \node(dynq1) at (3,-5.5) {$\vbar'_n(\bar{\Theta}_n,\Ybar_n)$};
      \draw[->] (in-EP)--(q1-EP) node[below=0.1cm,pos=0.5,color=blue!50] {\small quantization}; 
      \draw[->, thick,dotted] (va-EP)--(vaq1-EP) ;
      \draw[->, thick,dotted] (fi-EP)--(fiq1-EP);
	\draw[->, thick,dotted] (dyn)--(dynq1);
	\draw[->, thick,dotted]  (q1-EP)--(vaq1);
	\draw[->, thick,dotted]  (vaq1-EP)--(fiq1);
	\draw[->, thick,dotted] (fiq1-EP)--(dynq1);
      \node(fiq2) at (7,-3.5) {$(\hat{\Theta}_n,\Yhat_n)$};
       \node(fiq2-EP) at (7, -4) {$\Gamma_n$, $\hat{R}'_n$};
      \node(dynq2) at (7,-5.5) {$\vhat'_n(\hat{\Theta}_n,\Yhat_n)$};
      \draw[->] (fiq1-EP)--(fiq2-EP) node[below=0.15cm,pos=0.5,color=blue!50] {\small quantization};
      \draw[->, thick,dotted] (dynq1)--(dynq2) ;
      \draw[->, thick,dotted] (fiq2-EP)--(dynq2);
      \draw[->, thick,dotted] (dynq1)--(dynq2);
    \end{tikzpicture}}
  \end{center}
  \caption{Two-step approximation of the value functions}
  \label{fig:approx}
\end{figure}
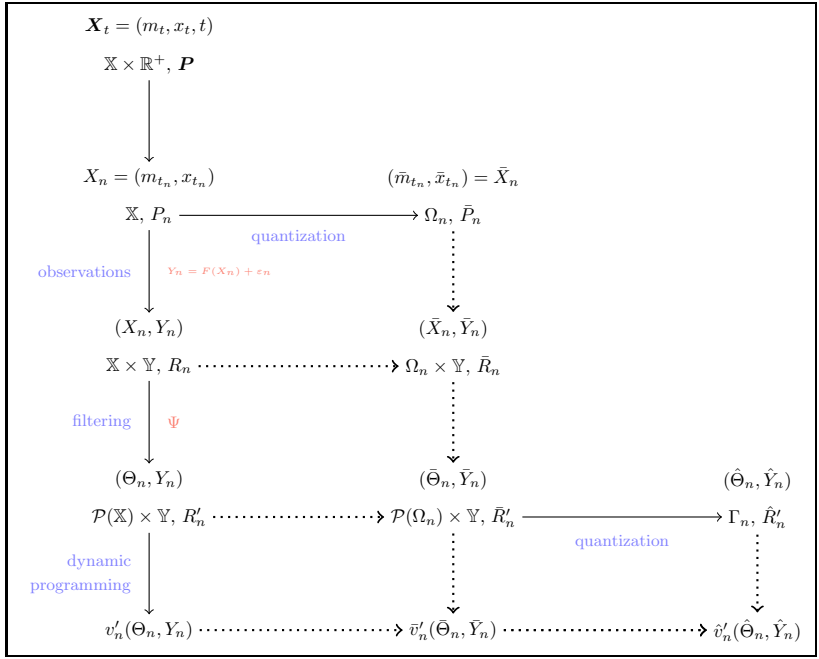

To build our approximation, we start from the dynamic programming equations from Theorem~\ref{th:DP'} and propose a two-step discretization of operators $R'_n$, $0~\leq~n~\leq~N~-~1$. Our global approach and the relationships between the different Markov chains we introduce, together with their state space and kernels are summarized in Figure \ref{fig:approx}. 
The left column corresponds to the construction presented in Section \ref{Framework} from the original continuous-time PDMP to the fully observed dynamic programming equations.

The first step in the middle column corresponds to a time-dependent discretization of the state space of the Markov chain~$(X_n)$. We obtain a finite state space Markov chain $(\bar{X}_n)$ that we plug into the observation equation~(\ref{Obs}) and filter operator to obtain Markov chains $(\bar{X}_n,\bar{Y}_n)$ with kernels $\bar{R}_n$ and $(\bar{\Theta}_n,\bar{Y}_n)$ with kernels $\bar{R}'_n$. Finally we replace $R'_n$ by $\bar{R}'_n$ in the dynamic programming equations to obtain the first sequence of approximate value functions. Note that by doing so, $\bar{Y}_n$ does not correspond to a discretization of the observations $Y_n$ and $\bar{\Theta}_n$ is not the filter of $\bar{X}_n$ given the observations ${Y}_n$. By this procedure, we start from a finite state space Markov chain $(\bar X_n)$ and obtain a simulatable filter $\bar{\Theta}_n$ that is still measure-valued but can be identified to finite-dimensional vectors. One more approximation is still required to obtain a finite state-space Markov chain.

The second step in the right column consists in the joint discretization of the Markov chain $(\bar{\Theta}_n,\bar{Y}_n)$. We obtain a finite state space Markov chain $(\hat{\Theta}_n,\hat{Y}_n)$ with kernel $\hat{R}'_n$. Again, we plug this new kernel into the dynamic programming equations. As the Markov chain $(\bar{\Theta}_n,\bar{Y}_n)$ has a finite state space, integrating with respect to $\hat{R}'_n$ simply corresponds to computing weighted sums. Hence the dynamic programming equations are now fully solvable numerically. This leads both to a numerically tractable approximation of the original value function and to a candidate $\epsilon$-optimal strategy.
Both steps are based on discretization by optimal quantization of the vector-valued Markov chains.

In this section,  we first introduce some notation and assumptions. Then we briefly recall the optimal quantization procedure and its main properties, proceed to construct the first discretization and state the convergence result of the approximate value functions to the original ones, and then construct the second discretization and state the convergence result of the approximate value functions to those from the previous part. Finally, we explain how a computable stopping strategy can be derived from the approximations of the value functions. The proofs of convergence are postponed to Section \ref{sec:proofs}.
%
\subsection{Notation and assumptions.}
\label{sec:notation}
In this section we introduce the function spaces we will be working with, the distance we consider for measure spaces and the main assumptions on our model.
Let $E$ be a Borel subset of $\X$.
%
\begin{defnt}
Let $BL(E)$ be the set of Borelian functions from $E \times \Y$ onto $\R$ for which there exist finite constants $\|\varphi\|_E$ and $[\varphi]_E$ such that for all $(m,x,y)$ and $(m,x',y)$ in $E\times\Y$, one has 
\begin{align*}
|\varphi(m,x,y)|&\leq \|\varphi\|_E,\\
|\varphi(m,x,y)-\varphi(m,x',y)|&\leq [\varphi]_E (|x-x'|).
\end{align*}
Denote also the unit ball of $BL(E)$ by
\begin{eqnarray*}
  BL_1(E) & = & \big\{ \varphi \in BL(E) : \|\varphi\|_E \leq 1 \ ; \ [\varphi]_E \leq 1  \big\}.
\end{eqnarray*}
\end{defnt}
%
\begin{defnt}
For $\theta$ and $\theta'$ two probability measures in $\mathcal{P}(E)$, define the distance $d_E(\theta, \theta')$ by
\begin{align*}
d_E(\theta, \theta')=
&\sup_{\varphi \in BL_1(E)}\sup_{y\in\Y} \left|\int\!\!\! \varphi(m,x,y) d\theta(m,x) -\!\!\int\!\!\! \varphi(m,x,y) d\theta' (m,x)\right|.
\end{align*}
\end{defnt}
%
In particular, if $E=\{e_1,\ldots, e_\ell\}$ is a finite set, this distance correspond to the $L_1$ distance on $\mathcal{P}(E)$:
$$d_E(\theta, \theta')=\|\theta- \theta'\|_{\mathcal{P}(E),1}=\sum_{i=1}^{\ell} \left| \theta(e_i)- \theta'(e_i) \right|.$$
%
\begin{defnt}
Let $BLP(E)$ be the set of Borelian functions from $\mathcal{P}(E) \times \Y$ onto $\R$ for which there exist finite constants $\|\varphi\|_{E,\mathcal{P}}$ and $[\varphi]_{E,\mathcal{P}}$ such that for all $(\theta,y)$ and $(\theta',y)$ in $\mathcal{P}(E)\times\Y$ one has
\begin{align*}
|\varphi(\theta,y)|&\leq \|\varphi\|_{E,\mathcal{P}},\\
|\varphi(\theta,y)-\varphi(\theta',y)|&\leq [\varphi]_{E,\mathcal{P}} d_E(\theta,\theta').
\end{align*}
\end{defnt}
%
In the sequel when $E=\X$ we will drop the index $E$:
for any function $\varphi$ in $BL(\X)$, $\|\varphi\|_{\X}=\|\varphi\|$ and $[\varphi]_\X=[\varphi]$, 
for any function $\varphi$ in $BLP(\X)$, $\|\varphi\|_{\X,\mathcal{P}}=\|\varphi\|_{\mathcal{P}}$ and $[\varphi]_{\X,\mathcal{P}}=[\varphi]_{\mathcal{P}}$, 
and for any probability measures $\theta$ and $\theta'$ in $\mathcal{P}(\X)$, $d_\X(\theta, \theta')=d(\theta, \theta')$.\\

We make the following regularity assumptions on the driving parameters of our processes. \\
\textbf{(H1)} There exist a function $L$ and a constant $L_{\!f}$ such that $\forall (x,{x'}) \in \K^2$ and $y\in\Y$ one has
\begin{align*}
\big|f(y-F(x))-f(y-F({x'}))\big|&\leq L(y)|x-x'|,\\
\int_\Y L(y) dy &\leq L_{\!f}<+\infty.
\end{align*}
%
\textbf{(H2)} There exist  positive real constants $\underline{f}$, $\overline{f}$ and $B_f$ such that for all $(x,y)\in\K\times\Y$,
\begin{align*}
0< \underline{f} \leq f(y-F(x))&\leq \overline{f}<+\infty,\\
\int_\Y\sup_{x\in\K}f(y'-F(x))dy'&\leq B_f<+\infty.
\end{align*}
%
\noindent \textbf{(H3)} For all $0\leq i \leq d$, there exists a positive constant $[\Phi_i]$ such that for all $0<t\leq N\delta$, one has
$$\big|\Phi_i({x},t)-\Phi_i({x'},t)\big|\leq  [\Phi_i] |x-{x'}|, \quad \forall (x,{x'}) \in \K^2.$$

\paragraph*{Examples} In the examples of Section \ref{sec:examples}, all the flows are clearly Lipschitz-continuous in $x$ for $0<t\leq N\delta$, so that \textbf{(H3)} holds.
For the identity link function $F$ and truncated Gaussian noise, one has
\begin{align*}
f(y-F(x))=\frac{1}{p\sigma\sqrt{2\pi}}e^{-\frac{(y-x)^2}{2\sigma^2}}\1_{[-s,s]}(y-x),
\end{align*}
with $p=\P(-s\leq Z\leq s)$ for a centered Gaussian random variable $Z$ with variance $\sigma^2$. Say the state space $\K$ is an interval of the form $[-S,S]$. Then $\Y=[-S-s, S+s]$ and Assumptions 
\textbf{(H1)} and \textbf{(H2)} hold with 
$L_{\!f}\leq 2sS(S+s)(p\sigma^3\sqrt{2\pi})^{-1}$, 
$\overline{f}\leq (p\sigma\sqrt{2\pi})^{-1}$, $B_f=2(s+S)\overline{f}$ and 
$\underline{f}\geq (p\sigma\sqrt{2\pi})^{-1}e^{-s^2/2\sigma^2}$.
Similarly, for the inverse link function and a state space $\K=[S_1,S_2]$ for some positive $S_1$, $S_2$, \textbf{(H2)} holds with the same constants $\underline{f}$ and $\overline{f}$, $B_f=(2s+S_2-S_1)\overline{f}$, and 
\textbf{(H1)} holds with 
$L_{\!f}\leq s(2s+S_2-S_1)(S_1^2p\sigma^3\sqrt{2\pi})^{-1}$.
%
\subsection{Quantization}
\label{sec:quantize}
%
We recall that for an $\mathbb{R}^{q}$-valued random variable $Z$ such that $\E[|Z|^2] < \infty$ and $\ell$ a fixed integer, the optimal $L_{2}$-quantization of the random variable $Z$ consists in finding the best possible $L_{2}$-approximation of $Z$
by a random variable $\widehat{Z}$ taking at most $\ell$ values in $\mathbb{R}^{q}$,
which will be denoted by $\Gamma_\ell=\{z^{1},\ldots,z^{\ell}\}$.
The asymptotic properties of the $L_{2}$-quantization are given by Zador's theorem (see, e.g., \cite[Theorem 3]{bally05}), recalled below, which uses the notation $p_{\Gamma}(z)$ for the closest neighbor projection of $z\in\mathbb{R}^{q}$ on a grid $\Gamma\subset\mathbb{R}^{q}$.
\begin{theor}\label{th-quantization}
Let $Z$ be an $\mathbb{R}^{q}$-valued random variable, and suppose that for some $\epsilon>0$ we have $\mathbb{E}[|Z|^{2+\epsilon}]~<~+~\infty$. Then, as $\ell$ tends to infinity, one has
\begin{align*}
 \min_{|\Gamma|\leq \ell} \E[ |Z-p_{\Gamma}(Z)|^{2}] = O(\ell^{-2/q} ).
\end{align*}
\end{theor}
There exist algorithms that can numerically find, for a fixed $\ell$, the quantization of~$Z$ (or, equivalently, the grid $\{z^{1}_{\ell},\ldots,z^{\ell}_{\ell}\}$ attaining the minimum in Theorem \ref{th-quantization} above and its distribution) as soon as $Z$ is simulatable \cite{pages98,pages04}. 
Roughly speaking, such a grid will have more points in the areas of high density of $Z$ and fewer points in the areas of low density of $X$. 
Replacing $Z$ by $\widehat{Z}$ turns integrals into finite sums and makes numerical computations possible, with easy derivation of error bounds for Lipschitz functionals of the random variable thanks to
Theorem \ref{th-quantization}.

Optimal quantization can also readily be extended to (discrete-time) Markov chains  \cite{pages04b}. One thus retrieves a quantization grid at each time step and the transition matrices between two consecutive grids. 
%
\subsection{First discretization}
\label{sec:first}
%
We propose a time-dependent discretization of the state space $\mathbb{X}$ based on the optimal quadratic quantization of the discrete time Markov chain $({X}_{n})_{n\geq 0}$. Let $\Omega_n$, $n\geq 0$ be a sequence of optimal quantization grids for $({X}_{n})_{n\geq 0}$. The cardinality of  $\Omega_n$ is denoted by $\ell_n$ and $\Omega_n=\{\omega_n^1,\dots,\omega_n^{\ell_n}\}$. Let $\bar{X}_n$ be the nearest-neighbor projection of ${X}_{n}$ onto $\Omega_n$: $\bar{X}_n=p_{\Omega_n}(X_n)$, and set
\begin{align*}
\bar P_n(\omega_{k+1}^j|\omega_n^i)=\bar p_{n,i,j}
&=\P(\bar{X}_{n+1}=\omega_{n+1}^j|\bar{X}_{n}=\omega_{n}^i).
\end{align*}
As the mode component is already discrete, we will assume in the sequel that the projection preserves the mode, i.e. if ${X}_{n}=(m,x)$ then $\bar{ X}_{n}=(m,\bar x)$.
To simplify notation, 
for any function $\varphi$ in $BL(\Omega_n)$, we will denote $\|\varphi\|_{\Omega_n}=\|\varphi\|_n$ and $[\varphi]_{\Omega_n}=[\varphi]_n$.

To define an approximation for the kernels $R'_n$, we replace in the definition of $R'_n$ the quantities related to the Markov chain $(X_n)$ by those related to $(\bar X_n)$. Namely, we define:\\
$\bullet$  a family of Markov kernels $\bar{R}_n$ from $\Omega_n\times\Y$ onto $\Omega_{n+1}\times\Y$:
\begin{align*}
\bar{R}_n(\omega_{n+1}^j,dy'|\omega_n^i,y)
&=\bar{R}_n(\omega_{n+1}^j,dy'|\omega_n^i)\\
&=\bar p_{n,i,j}f(y'-F(\omega_{n+1}^j))dy',
\end{align*}
$\bullet$  a family of operators $\bar{\Psi}_n$ from $\mathcal{P}(\Omega_n)\times \mathbb{R}$ onto $\mathcal{P}(\Omega_{n+1})$:
\begin{align*}
{\Psibar_n(\theta,y')(\omega_{n+1}^j)}=\frac{\sum_{i=1}^{\ell_{k}}\bar p_{k,i,j}\fF_{y'}(\omega_{k+1}^j){\theta}(\omega_{k}^i)}{\sum_{j'=1}^{\ell_{k+1}}\sum_{i=1}^{\ell_{k}}\bar p_{k,i,j'}\fF_{y'}(\omega_{k+1}^{j'}){\theta}(\omega_{k}^i)}.
\end{align*}
%
From these two ingredients, we construct a new family of Markov kernels $\bar{R}_n'$ from $\mathcal{P}(\Omega_n)\times\mathbb{Y}$ onto $\mathcal{P}(\Omega_{n+1})\times\mathbb{Y}$ by setting, for all Borelian subsets $\bar P\subset\mathbb{P}(\Omega_n)$ and $\bar C\subset \Y$,
\begin{align*}
\bar{R}'_n(\bar{P}\times\bar{C}|\thetabar,\bar{y})=\bar{R}'_n(\bar{P}\times\bar{C}|\thetabar)
  &=\!\!\!\!\!\!\sum_{\substack{1\leq i\leq \ell_{n}\\1\leq j\leq \ell_{n+1}}}\!\!\!\!\!\!\int_{\bar{C}}\!\!\! \1_{\bar{P}}(\bar{\Psi}_n(\thetabar,\bar{y}'))\bar p_{n,i,j}f(y'-F(\omega_{n+1}^j))dy'\thetabar(\omega_n^i).
  \end{align*}
From the family of Markov kernels $(\bar R'_n)$, one can construct a Markov chain $(\bar\Theta_n,\bar Y_n)_{n\geq 0}$ by setting
\begin{align*}
\P(\bar\Theta_0=\delta_{(0,x_0)},\bar Y_0=\bar y)&=1\\
\P((\bar\Theta_{n+1},\bar Y_{n+1})\in \bar P\times \bar C\ | \bar\Theta_n=\bar\theta,\bar Y_n=\bar y)&=\bar R'_n(\bar P\times \bar C|\bar\theta).
\end{align*}
Note that \\
$\bullet$ $(\bar Y_n)$ does not have the same dynamics as the original observations $(Y_n)$, nor does it correspond to a function of these observations;\\
$\bullet$ one has $\bar \Theta_n(\omega_n^i)=\P(\bar X_n =\omega_n^i|\bar Y_0,\ldots \bar Y_n)$ which thus corresponds to the filter of $(\bar X_n)$ given the $(\bar Y_n)$, but does not correspond to the filter of $(\bar X_n)$ given the original observations $(Y_n)$, or to a function of the original filter $(\Theta_n)$ nor of observations $(Y_n)$;\\
$\bullet$ unlike the recursion for the filter $\Theta_n$, the recursion for $\bar\Theta_n$ is numerically tractable, thus $(\bar\Theta_n)$ is simulatable.

The random filter $\bar\Theta_n$ is characterized by the random weights $\Thetabar_n^i=\Thetabar_n(\omega_n^i)$, for $i=1,\dots,\ell_n$ and can be identified with a random vector valued in the $\ell_n$-simplex
in $\R^{\ell_n}$ of dimension $\ell_n-1$. 
This identification will be in force throughout this paper.

 Finally,  we define the main quantities of interest for this section, namely the approximate value functions $\vbar_n'$ from $\mathcal{P}(\Omega_n)\times \Y$ onto $\mathbb{R}$ as
    \begin{eqnarray*}
      \vbar_N'(\thetabar,y')&=& \min_{a\in\mathcal{M}} \ C'(\thetabar, y',a) \\
      \vbar_n'(\thetabar,y')&=&\min\left\{\min_{a\in\mathcal{M}_+} C'(\thetabar,y',a); c'(\thetabar,y')+\Rbar'_n\vbar'_{n+1}(\thetabar,y') \right\}.
    \end{eqnarray*}
Then we have the following convergence.
\begin{theor}\label{lem:v_n} 
Under assumptions \textbf{(H1-3)}, for all distributions $\theta$ in $\mathcal{P}(\X)$, $\thetabar$ in $\mathcal{P}(\Omega_N)$ and all $(y,\ybar)\in\mathbb{Y}^2$, one has
  \begin{align*}
    |v'_N(\theta, y) - \vbar_N'(\thetabar,\ybar)|     \leq \bar B \ d(\theta,\thetabar) .
  \end{align*}
For all $0\leq n\leq N-1$, all distributions $\theta$ in $\mathcal{P}(\X)$, $\thetabar$ in $\mathcal{P}(\Omega_n)$ and all $(y,\ybar)\in\mathbb{Y}^2$, one has
\begin{align*}
|v'_n(\theta, y) - \vbar_n'(\thetabar,\ybar)|
&\leq [v'_n]_{\mathcal{P}} d(\theta,\thetabar) + \|\vbar'_{n+1} -v'_{n+1}\|_{{n+1}}+\big( L_{\!f}\|v'_{n+1}\|_{{\mathcal{P}}}+\overline{f}\ \fFbar [v'_{n+1}]_{{\mathcal{P}}}\big)\\
&\quad\times(\E[|\bar{X}_{n+1}-{X}_{n+1}|]+[\Phi]^2\E[|X_n-\bar{X}_n|]).
\end{align*}
\end{theor}
%
The proof of this theorem is given in Section \ref{sec:proof-theo}, and the constants $\bar B$ and $\overline{f\!F}$ are given in Table~\ref{tab:cstes}.
In particular, for the optimal performance one has
  \begin{align*}
|v'_0(\delta_{(0,x_0)},y_0)- \vbar_0'(\delta_{(0,x_0)},y_0)|
&\leq  \sum_{n=0}^{N-1}a_n\E[|\bar{X}_{n}-{X}_{n}|],
  \end{align*}     
with 
\begin{eqnarray*}
a_n&=&L_{\!f}(\|v'_{n}\|_{\mathcal{P}}+[\Phi]^2\|v'_{n+1}\|_{\mathcal{P}})+\overline{f}\ \fFbar([v'_{n}]_{\mathcal{P}}+[\Phi]^2[v'_{n+1}]_{\mathcal{P}})\\
&\leq &\bar B(1+[\Phi]^2)\Big(L_{\!f}(N-n+1)+\overline{f}\ \fFbar\sum_{j=0}^{N-n}\big((1+[\Phi]^2L_{\!f})(N-n-j)+1\big)\!(\fFbarp\overline{f})^j\Big).
\end{eqnarray*}
%
\subsection{Second discretization}
\label{sec:second}
%
The value functions $\bar v'_n$ are not directly numerically computable as they involve integration by operators $\bar R'_n$ on the continuous space $\mathcal{P}(\Omega_n)\times\Y$. A second discretization is thus needed.
To do so, we use optimal quantization again. Following the generic direction given in \cite{pham05}, we discretize jointly $(\bar \Theta_n,\bar Y_n)$. Note that $(\bar \Theta_n,\bar Y_n)$ is easy to simulate because the recursive construction of $\bar\Theta_n$ only involves finite weighted sums. This approach would not have been possible on the chain $(\Theta_n,Y_n)$ as $\Theta_n$ cannot be simulated exactly.

Our second discretization step thus now consists in replacing the Markov chain $(\bar\Theta_n,\Ybar_n)$ by its optimal quantization approximation $(\hat\Theta_n,\Yhat_n)$.
By construction, $(\hat\Theta_n,\Yhat_n)$ takes a finite number of values on a grid $\Gamma_n$ of size $N_n$:
$\Gamma_n=\{\gamma_n^i\}_{1\leq i\leq N_n} = \{(\pi_n^1,y_n^1),\dots,(\pi_n^{N_n},y_n^{N_n}) \}.$
To simplify notation, 
for any measure $\theta$ in $\mathcal{P}(\Omega_n)$, we will simply denote $\|\theta\|_{\mathcal{P}(\Omega_n),1}=\|\theta\|_{n,1}$.

We now set
\begin{eqnarray*}
  \Rhat'_n(\gamma_{n+1}^j|\gamma_n^i)
  &=&\P((\hat\Theta_{n+1},\Yhat_{n+1})=\gamma_{n+1}^j|(\hat\Theta_n,\Yhat_n)=\gamma_{n}^i).
\end{eqnarray*}
With these new transition kernels, we define the approximate value functions $\vhat_n'$ from $\Gamma_n$ onto $\R$ as
  \begin{eqnarray*}
    \vhat_N'(\gamma)&=& \min_{a\in\mathcal{M}} \ C'(\gamma,a), \\
    \vhat_n'(\gamma)&=&\min\left\{\min_{a\in\mathcal{M}_+} C'(\gamma,a); c'(\gamma)+\Rhat'_n\vhat'_{n+1}(\gamma) \right\}.
  \end{eqnarray*}
%
%
These functions can be numerically computed on the grids $\Gamma_n$.
%
\begin{theor} \label{lem:vbark-vhatk} 
Under assumptions \textbf{(H1-3)}, for all distribution $\thetabar_N$ in $\mathcal{P}(\Omega_N)$, $\ybar \in \Y$ and $(\thetahat_N,\yhat) \in \Gamma_N$, one has
  \begin{eqnarray*}
    |\vhat'_N(\thetahat_N, \yhat) - \vbar_N'(\thetabar_N,\ybar)|  \leq  \bar B \|\thetabar_N-\thetahat_N\|_{{N,1}}.
  \end{eqnarray*}
For $0\leq n\leq N-1$, for all distribution $\thetabar_n$ in $\mathcal{P}(\Omega_n)$, $\ybar \in \Y$ and $(\thetahat_n,\yhat)\in \Gamma_n$, one has
    \begin{align*}
|\vhat'_n(\thetahat_n, \yhat) - \vbar_n'(\thetabar_n,\ybar)| & \leq [\vbar'_n]_{n} \|\thetahat_n-\thetabar_n\|_{n,1} + \|\vhat'_{n+1}-\vbar'_{n+1}\|_{{n+1}}  + [\vbar'_{n+1}]_{n+1}\E[\|\Thetahat_{n+1}-\Thetabar_{n+1}\|_{n+1,1}] \\
      & + (\|\vbar'_{n+1}\|_{n+1}+\overline{f}\fFtilde [\vbar'_{n+1}]_{n+1}) \E[\|\Thetahat_{n}-\Thetabar_{n}\|_{n,1}],
    \end{align*}
 \end{theor}

The proof of this theorem is given in Section \ref{sec:proof-theo}, and the constant $\fFtilde$ is given in Table~\ref{tab:cstes}.
 In particular, we obtain
 \begin{align*}
   |\vhat'_0(\delta_{(0,x_0)},y_0) - \vbar_0'(\delta_{(0,x_0)},y_0)|    & \leq   \sum_{n=0}^{N} b_n \E[\dn{\Thetahat_{n}-\Thetabar_{n}}],
 \end{align*}
 with
 \begin{align*}
   b_n &= \barnlip{\vbar'_{n}} + \barnpsup{\vbar'_{n+1}}+\overline{f}\fFtilde \barnplip{\vbar'_{n+1}} \\
    &\leq 2 \bar B \sum_{j=0}^{N-n}(N-n+1-j)(\overline{f} \fFtilde)^j .
\end{align*}
%
Suppose that all quantization grids $\Omega_n$ have the same number of points $N_\Omega$ and all quantization grids $\Gamma_n$ also have the same number of points $N_\Gamma$. Then, Theorem \ref{th-quantization} yields $\E[|X_n-\bar X_n|]=O(N_\Omega^{-1})$ and $\E[|\Yhat_{n}-\Ybar_{n}|+\dn{\Thetahat_{n}-\Thetabar_{n}}]=O(N_\Gamma^{-1/N_\Omega})$. Thus one has
 \begin{align*}
 |v_0'(\delta_{(0,x_0)},y_0)-\vhat'_0(\delta_{(0,x_0)},y_0)|
= O(N_\Omega^{-1}+N_\Gamma^{-1/N_\Omega}).
 \end{align*}
This rate of convergence is very slow, which is not surprising given that one had to discretize infinite dimension measure-valued random variables. This is the well known curse of dimensionality one is faced with when dealing with partial observations.
%
\subsection{Construction of a stopping strategy.}
\label{sec:stop}
%
We can now construct a computable stopping strategy using the fully discretized value function.
Suppose that the process starts from point $\xi_0=(0,x_0,y_0)$ and observations $y_0,\ldots,y_n$ are available at time $n$. One cannot compute the filter $\P_{\xi_0}(X_n\in \cdot | (Y_0,\ldots,Y_n)=(y_0,\ldots,y_n))$ because of the continuous integrals in the definition of the Bayes operator from Proposition \ref{def:psi}. However, one can recursively compute an approximate filter as follows.
\begin{align*}
\bar\theta_0&=\delta_{(0,x_0)},\quad\bar\theta_k=\bar\Psi_{k-1}(\bar\theta_{k-1},y_k),\ 1\leq k\leq n.
\end{align*}
By construction, $\bar\theta_k$ belongs to $\mathcal{P}(\Omega_k)$ for all $k$. Then this approximate filter can be projected onto the quantization grids $(\Gamma_k)_{0\leq k\leq n}$: 
\begin{align*}
(\hat\theta_k,\hat y_k)&=p_{\Gamma_k}(\bar\theta_k,y_k),
\end{align*}
for all $1\leq k\leq n$. Finally, the values of $\hat v'_k(\hat\theta_k,\hat y_k)$ are available for all $0\leq k\leq n$.\\

Now we define two sequences of function $(r_n)_{0\leq n\leq N}$ and $(a_n)_{0\leq n\leq N}$ as\\
$\bullet$ for $0\leq n\leq N-1$, $r_n : \mathcal{P}(\Omega_n)\times\Y\to \{0,1\}$ and $a_n : \mathcal{P}(\Omega_n)\times\Y\to \mathcal{M}_+$ are such that 
\begin{align*}
r_n(\bar \theta,\bar y)&=\1_{\min\limits_{a\in\mathcal{M}_+}C'(p_{\Gamma_n}(\bar \theta,\bar y),a)<c'(p_{\Gamma_n}(\bar \theta,\bar y))+\hat R'_n\hat v_{n+1}(p_{\Gamma_n}(\bar \theta,\bar y))},\\
a_n(\bar \theta,\bar y)&=\arg\min_{a\in\mathcal{M}_+}C'(p_{\Gamma_n}(\bar \theta,\bar y),a)\1_{(r_n(\bar \theta,\bar y)>0)}.
\end{align*}
$\bullet$ $r_N : \mathcal{P}(\Omega_N)\times\Y\to \{0,1\}$ and $a_N : \mathcal{P}(\Omega_N)\times\Y\to \mathcal{M}$ are such that 
\begin{align*}
r_N(\bar \theta,\bar y)&=\1_{(\arg\min_{a\in\mathcal{M}}C'(p_{\Gamma_N}(\bar \theta,\bar y),a)>0)},\\
a_N(\bar \theta,\bar y)&=\arg\min_{a\in\mathcal{M}}C'(p_{\Gamma_N}(\bar \theta,\bar y),a).
\end{align*}
%
Thus $r_n$ is a stopping indicator depending on which term won the minimization in the dynamic programming, and $a_n$ corresponds to the mode to be selected after the jump.
\begin{figure}[t]
\tikzstyle{decision} = [diamond, draw, 
    text width=4.5em, text badly centered, node distance=3cm, inner sep=0pt]
\tikzstyle{block} = [rectangle, draw, text width=15em, 
text centered, rounded corners, minimum height=4em]
\tikzstyle{line} = [draw, ->]
\tikzstyle{cloud} = [draw, ellipse,node distance=3cm,
    minimum height=2em]

\begin{center}    
\framebox{\begin{tikzpicture}[node distance = 2cm, auto, scale=1, every node/.style={scale=0.65}]
    \node [block] (init) {$n\leftarrow 0$\\ $y\leftarrow y_0$\\ $\bar\theta\leftarrow\delta_{0(,x_0)}$\\$r\leftarrow r_0(\bar\theta,y)$};
    \node [cloud, right of=init, xshift=3cm] (obs0) {Observation $y_0$};
    \node [decision, below of=init] (decide) {$r=1$ ?};
    \node [block, right of=decide, xshift=3cm] (stop) {Stop at time $n$\\Choose decision $a=a_n(\bar\theta,y)$};
    \node [decision, below of=decide, node distance=3cm] (horizon) {$n=N$ ?};
    \node [block, right of=horizon, xshift=3cm] (stopN) {Choose decision $a=0$};
    \node [block, below of=horizon, node distance=3cm] (update) {$n\leftarrow n+1$\\ $y\leftarrow y_{n}$\\ $\bar\theta\leftarrow\Psi_{n-1}(\bar\theta,y)$\\$r\leftarrow r_n(\bar\theta,y)$};
    \node [cloud, right of=update, xshift=3cm] (obsn) {Observation $y_n$};
    \path [line] (init) -- (decide);
    \path [line] (decide) -- node {yes}(stop);
    \path [line] (decide) -- node {no}(horizon);
    \path [line] (horizon) -- node {yes}(stopN);
    \path [line] (horizon) -- node {no}(update);   
    \path [line] (update) -- +(-4,0) |- (decide);
    \path [line,dashed] (obs0) -- (init);
    \path [line,dashed] (obsn) -- (update);
\end{tikzpicture}}
\end{center}
\caption{Computable stopping strategy}
\label{fig:strategy}
\end{figure}
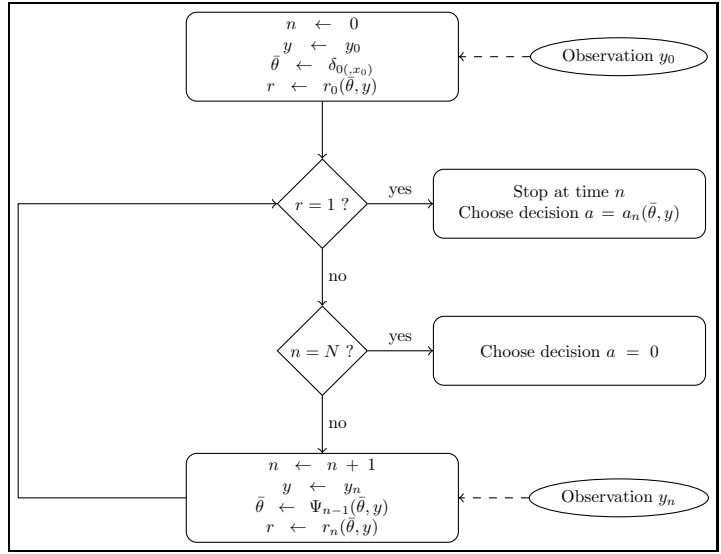

Our candidate stopping strategy is the following, as illustrated on Figure \ref{fig:strategy}:\\
Compute $r_0(\delta_{(0,x_0)},y_0)$\\
$\bullet$ if $r_0(\delta_{(0,x_0)},y_0)=1$, stop at time $0$ and select decision $a_0(\delta_{(0,x_0)},y_0)$\\
$\bullet$ otherwise compute $\bar\theta_1=\bar\Psi_{0}(\delta_{(0,x_0)},y_1)$ and $r_1(\bar\theta_1,y_1)$\\
\phantom{\quad} $-$ if $r_1(\bar\theta_1,y_1)=1$, stop at time $1$ and select decision $a_1(\bar\theta_1,y_1)$\\
\phantom{\quad} $-$ otherwise compute $\bar\theta_2=\bar\Psi_{1}(\bar\theta_1,y_2)$ and $r_2(\bar\theta_2,y_2)$ and so on until time $N$\\
\phantom{\quad\quad} $\cdot$ if the process wasn't stopped before $N$, compute $\bar\theta_N=\bar\Psi_{N-1}(\bar\theta_{N-1},y_N)$ and $r_N(\bar\theta_N,y_N)$\\
\phantom{\quad\quad} $\cdot$ if $r_N(\bar\theta_N,y_N)=1$, stop at time $N$ and select decision $a_N(\bar\theta_N,y_N)$\\
\phantom{\quad\quad} $\cdot$ otherwise select decision $a=0$.

Note that all quantities can be computed numerically and that this strategy is non-anticipative. However, the sequence $(\bar\theta_n,y_n)$ is not a realization of the Markov chain $(\bar \Theta_n,\bar Y_n)$ nor of $(\Theta_n,Y_n)$. Therefore assessing theoretically the performance of this strategy is an open question that will be the subject of future works. Its numerical performance is assessed in Section~\ref{Simulations}.
%
\section{Convergence of the approximations}
\label{sec:proofs}
%
This section is dedicated to the proof of Theorems \ref{lem:v_n} and \ref{lem:vbark-vhatk}. 
Some useful constants for the sequel are given in Table~\ref{tab:cstes}.
\begin{table}[h!]
\centering
\begin{tabular}{ll}
\hline
$ [\Phi] = 1 \vee \max_{0\leq i\leq d}\{\phiilip\}$   &$\bar G = \alpha \vee \max_{(m,a)\in\mathcal{M}_+^2}\{\gamma_{ma}\}$\\
$\bar B =\bar G + \delta \beta$ 			     &$\overline{f\!F}= (B_{\!f}+L_{\!f}){\underline{f}}^{-1} + L_{\!f} \overline{f}\underline{f}^{-2}$\\
$\fFtilde=  B_{\!f}\underline{f}^{-1}(1+ \overline{f}\underline{f}^{-1})$ &  $B(y')=\sup_{x\in\K}f(y'-F(x))$\\
\multicolumn{2}{l}{$\overline{f\!F}^+= B_{\!f}\underline{f}^{-1}+[\Phi]^2\overline{f\!F}$}\\
\hline
\end{tabular}
\caption{{ Lipschitz-regularity constants and bounds}} \label{tab:cstes}
\end{table} 
%
\subsection{Regularity of the Markov operators and value functions}
We start with regularity properties of the cost functions, kernels and operators involved in the discretizations. The first lemma is straightforward and its proof is omitted.
%
\begin{lemma}\label{lem:ga-BL}
  For all $a\in\mathcal{M}_+$, set $C_a(m,x,y)=C_a(m)=C(m,x,y,a)$, then
$C_a$,
    and $c$ are in $BL(\X)$ with
  \begin{equation*}
  \|C_a\|=\bar G,\quad[C_a]=0,\quad
\|c\|= \delta\beta,\quad [c]=0, 
\end{equation*}
 $C_a'$ and $c'$ are in $BLP(\X)\cap BLP(\Omega_k)$ for all $1\leq k \leq N$ with
  \begin{equation*}
    \begin{array}{ll}
  \|C_a'\|_{\mathcal{P}}=\|C_a'\|_{k}   \leq \bar G,&  \quad {[C_a']}_{\mathcal{P}}={[C_a']}_{k}   \leq \bar G, \\
  \|c'\|_{\mathcal{P}}= \|c'\|_{k}    \leq \delta\beta, &  \quad{[c']}_{\mathcal{P}}={[c']}_{k}  \leq \delta\beta.
  \end{array}
\end{equation*}
\end{lemma}
%
In the sequel, for any function $g:\X\times\Y\rightarrow\mathbb{R}$ we will still denote $P_ng$ the function $\X\times\Y$ onto $\mathbb{R}$ defined by 
$$P_ng(m,x,y)=\E[g(X_{n+1},y)|X_n=(m,x)].$$
%
\begin{lemma}\label{lem:Pflip} 
Under assumption \textbf{(H3)}, for any $g \in BL(\X)$, $P_kg$ is in $BL(\X)$ with 
\begin{align*}
\|P_kg\|\leq \|g\|,\quad [P_kg]\leq \philip^2 [g].
\end{align*}
\end{lemma}
%
\noindent\textbf{Proof}
As $P_k$ is a Markov kernel, one clearly has $\|P_kg\|\leq \|g\|$. On the one hand, if $m=i> 0$, for $y\in\Y$, one has
\begin{align*}
\left| P_kg(i,x,y)- P_kg(i,\bar{x},y)\right|
&= |g(i,\Phi_i(x,\delta),y)-g(i,\Phi_i(\bar x,\delta),y)|  \leq[g] [\Phi_i] |x-\bar x|.
\end{align*}
On the other hand, if $m=0$, one has
\begin{align*}
\lefteqn{\left|P_kg(0,x,y)-P_kg(0,\bar{x},y)\right|}\\
&\leq  e^{-\int_0^\delta\lambda(k\delta+s)ds } |g(0,\Phi_0(x,\delta),y)-g(0,\Phi_0(\bar x,\delta),y)| \\
&\ +   \sum_{i=1}^d \pi_i \int_0^\delta \lambda(k\delta+s)e^{-\int_0^s\lambda(k\delta+z)dz }\Big|g(i,\Phi_i(\Phi_0(x,s),\delta-s),y)-g(i,\Phi_i(\Phi_0(\bar x,s),\delta-s),y)\Big|ds\\
&\leq[g] \phiilip \phi0lip |x-\bar x|,
\end{align*}
and we conclude using the definition of $\philip$.  
\hfill$\Box$
%
 \begin{lemma}\label{lem:Psi-dist}
Under assumptions \textbf{(H1-3)}, for all $0\leq k\leq N-1$ and $\theta\in\mathcal{P}(\X)$, $\thetabar\in\mathcal{P}(\Omega_k)$, one has
\begin{align*}
   \int_\Y d(\Psi_k(\theta,y'),\Psi_k(\thetabar,y'))dy'\leq \fFbarp \ d(\theta,\thetabar).
\end{align*}
Under assumptions \textbf{(H2-3)}, for all $(\theta,\thetabar) \in \mathcal{P}(\Omega_{k})^2$  one has
 \begin{eqnarray*}
\int_\Y\dkp{\Psibar_{k}(\theta,y')-\Psibar_{k}(\thetabar,y')}dy'  \leq \ \fFtilde \ \dk{\theta-\thetabar}.
 \end{eqnarray*}
\end{lemma}
%
\noindent\textbf{Proof} 
Set $g\in BL_1(\X)$. For $(y,y')\in\Y^2$, we set  $ f\!F_{y'}[g]:(m,x,y)\mapsto g(m,x,y)f(y'-F(x))$. Thus, one has
\begin{align*}
\lefteqn{\int_\X g(\xi',y)\ d\Psi_k(\theta,y')(\xi')-\int_\X g(\xi',y)\ d\Psi_k(\thetabar,y')(\xi')}\\
&=\frac{\d\int_\X P_k f\!F_{y'}[g](\xi,y)d\theta(\xi)-\d\sum_{i=1}^{\ell_{k}}P_k f\!F_{y'}[g](\xi,y)\thetabar(\omega_k^i)}{\d\int_\X P_kf\!F_{y'}(\xi)d\theta(\xi)} \\
&+ \d\sum_{i=1}^{\ell_{k}}P_k f\!F_{y'}[g](\xi,y)\thetabar(\omega_k^i) \Bigg(\Big(\d\int_\X P_kf\!F_{y'}(\xi)d\theta(\xi)\Big)^{-1}\!\!\!\!\!-\Big(\d\sum_{i=1}^{\ell_{k}} P_{k}f\!F_{y'}(\omega_k^i)\thetabar(\omega_k^i)\Big)^{-1} \Bigg).
\end{align*}
By assumption \textbf{(H2)} all denominators are bounded from below by $\underline{f}$. 
Next, the function $ f\!F_{y'}[g]$ 
is clearly in $BL(\X)$ with $\|f\!F_{y'}[g]\|  \leq B(y')$ and ${[f\!F_{y'}[g]]} \leq L(y')+ B(y')$ as $g\in BL_1(\X)$. Hence, 
by Lemma \ref{lem:Pflip} we readily have $P_kf\!F_{y'}[g] \in BL(\X)$ 
so that for all $g\in BL_1(\X)$ we have
\begin{align*}
\lefteqn{\int_\X g(\xi',y)\ d\Psi_k(\theta,y')(\xi')-\int_\X g(\xi',y)\ d\Psi_k(\thetabar,y')(\xi')}\\
& \leq \big((B(y')\!+\!\philip^2 (B(y') \!+\! L(y'))\underline{f}^{-1}+B(y')\philip^2L(y')\underline{f}^{-2}\big)d(\theta,\thetabar),
\end{align*}
leading to the expected result after integrating on $y'$. The proof of the second statement follows the same lines and is omitted.
\hfill$\Box$
%
\begin{lemma}\label{Rk'-lip}
Let $g \in BLP(\X)$. For all $0\leq k\leq N-1$, under assumptions \textbf{(H1-3)}, $R_k'g$ is in $BLP(\X)$ with 
\begin{align*}
\|R_k'g\|_{\mathcal{P}}&\leq \|g\|_{\mathcal{P}},\\
   [R_k'g]_{\mathcal{P}}&\leq (1+[\Phi]^2L_{\!f})\|g\|_{\mathcal{P}} + \fFbarp \overline{f} [g]_{\mathcal{P}}.
\end{align*}
Let $0\leq k\leq N-1$ and $g$ in $BLP(\Omega_{k+1})$. Then, under assumptions \textbf{(H2-3)}, $\Rbar'_{k} g$ is in $BLP(\Omega_{k})$ with
 \begin{align*}
  \|\Rbar'_{k} g\|_{k}  & \leq \|g\|_{k+1}, \\
  [\Rbar'_{k} g ]_{k} &\leq   [g]_{k+1}\fFtilde\overline{f} +\|g\|_{k+1}.
 \end{align*}
\end{lemma}
%
\noindent \textbf{Proof}
 As we deal with Markov kernels, we clearly have $\|R_k'g\|_{\mathcal{P}}\leq \|g\|_{\mathcal{P}}$ and $ \|\Rbar'_{k} g\|_{k} \leq \|g\|_{k+1}$. 
First, let us prove that if $g \in BLP(\X)$, then $\bar g_\theta$ defined from $\X\times\Y$ onto $\R$ by
\begin{align*}
{\bar g_\theta(m,x,y)}
&=\int_\Y g(\Psi_k(\theta,y'),y')P_kf\!F_{y'}(m,x)dy'
\end{align*}
is in $BL(\X)$.
On the one hand, one clearly has $\|\bar g_\theta\|\leq \|g\|_{\mathcal{P}}$.
On the other hand, using similar computations as in Lemma \ref{lem:Pflip} for $f\!F_{y'}$, one gets
\begin{align*}
|\bar g_\theta(m,x,y)-\bar g_\theta(m,\bar x,y)| 
&\leq  \int_\Y |g(\Psi_k(\theta,y'),y')||P_{k}\fF_{y'}(m,x)-P_{k}\fF_{y'}( m,\bar x)|dy'\\
&\leq \|g\|_{\mathcal{P}}[\Phi]^2L_{\!f}|x-\bar x|.
\end{align*}
Now, getting back to $R'_k$,
By definition, $R_k'g(\theta,y)$ does not depend on $y$ and one has
\begin{align*}
R_k'g(\theta,y)
&=\int_\X \bar g_{\theta}(m,x)d\theta(m,x),
\end{align*}
which yields
\begin{align*}
|R_k'g(\theta,y)-R_k'g(\bar\theta,y)|
&\leq (\|\bar g_\theta\| +[\bar g_\theta])d(\theta,\bar \theta)+\int_\X |\bar g_\theta (m,x)- \bar g_{\bar\theta}(m,x)|d\bar\theta(m,x).
\end{align*}
Finally, one uses Lemma \ref{lem:Psi-dist} to obtain
\begin{align*}
|\bar g_\theta (m,x)- \bar g_{\bar\theta}(m,x)|
&\leq [g]_{\mathcal{P}}\int_\Y  d\left(\Psi_k(\theta,y')),\Psi_k(\bar\theta,y')\right)P_{k}\fF_{y'}(m,x)dy' \\
&\leq [g]_{\mathcal{P}}\overline{f} \ \fFbarp d(\theta,\thetabar).
\end{align*}
Combining all the results, one gets the expected bounds. To obtain the Lipshitz property on $BLP(\Omega_k)$ one applies the same decomposition and the second statement of Lemma \ref{lem:Psi-dist}.
\hfill$\Box$
%
\begin{lemma}\label{stability'} 
Under assumptions \textbf{(H1-3)}, for all  $0\leq n\leq N$,  $v'_{n} \in BLP(\X)$ with
\begin{eqnarray*}
  \|v'_{n}\|_{\mathcal{P}}   &\leq&(N-n+1)\bar B,\\
  {[v'_{n}]}_{\mathcal{P}}  &\leq &  \bar B\!\!\sum_{k=0}^{N-n}\!\!\big((1+[\Phi]^2L_{\!f})(N-n-k)+1\big)\!(\fFbarp\overline{f})^k\!\!.
 \end{eqnarray*}
Under assumptions {\bf (H2-3)}, for all $1\leq n\leq N$, $\vbar_n$ is in $BLP(\Omega_n)$ with
  \begin{align*}
    \|\vbar_{n}\|_{n}& \leq (N-n+1) \bar B, \\
    [\vbar_{n}]_n& \leq \bar B \sum_{k=0}^{N-n}\big(N-n-k+1\big)(\fFtilde\overline{f})^k.
  \end{align*}
\end{lemma}
%
\noindent \textbf{Proof}
We prove here the first part of the statement, proceeding  by backward induction on $k$ using Lemma~\ref{lem:ga-BL}: \\
$\bullet$ For $k=N$, let $(\theta,y)$ and $(\bar\theta,\bar y) \in \mathcal{P}(\X)\times \Y$
\begin{align*}
  |v'_N(\theta,y)-v'_N(\bar\theta,y)| &=\left|\min_{a\in\mathcal{M}} C'(\theta,y,a)-\min_{a\in\mathcal{M}} C'(\bar\theta,y,a)\right| \leq \bar B  \ d(\theta,\bar\theta),\\
|v'_N(\theta,y)|&=\big|\min_{a\in\mathcal{M}} C'(\theta,y,a)\big| \leq \max_{a\in\mathcal{M}_+} \left|C_a'(\theta,y,a)\right|\vee |c'(\theta,y)| \leq \bar B.
\end{align*}
$\bullet$ We now assume that the statement holds for a given $k+1$ in $\{1,\dots,N\}$. 
According to Lemma \ref{Rk'-lip}, as $v'_{k+1}$ is in $BLP(\X)$. 
From the definition of $v'_k$, one thus obtains
\begin{eqnarray*}
  \|v'_{k}\|_{\mathcal{P}}  &\leq&\max_{a\in\mathcal{M}_+} \|C_a'\|_{\mathcal{P}}\vee\big( \|c'\|_{\mathcal{P}}+ \|v'_{k+1}\|_{\mathcal{P}}\big)\\
  &\leq
  & \bar B  +\|v'_{k+1}\|_{\mathcal{P}}\leq  (N-k+1)\bar B.
\end{eqnarray*}
Regarding the Lipschitz property, one has
\begin{align*}
|v'_k(\theta,y)-v'_k(\thetabar,y)|&\leq\max_{a\in\mathcal{M}_+}\left| C_a'(\theta,y)- C_a'(\thetabar,y)\right|\vee\left\{|c'(\theta,y)-c'(\bar\theta,y)|+|R'_kv'_{k+1}(\theta,y) -R'_kv'_{k+1}(\thetabar,y)| \right\}\\
  &\leq \!\big(\!\bar B\! +\!\|v'_{k+1}\|_{\mathcal{P}}(1\!+\![\Phi]^2L_{\!f})\!+\![v'_{k+1}]_{\mathcal{P}} \fFbarp\overline{f}\big) d(\theta,\bar\theta),
 \end{align*} 
using Lemma \ref{Rk'-lip} again.
\hfill$\Box$
 %
\subsection{Comparison of the dynamic programming quantities}
We start with a technical lemma comparing the Markov chains $(X_n)$ and $(\bar X_n)$.
\begin{lemma}\label{lem:r-rbar integral}
Under assumptions \textbf{(H1)} and \textbf{(H3)}, for all $0\leq k\leq N-1$ and $1\leq i\leq \ell_k$, one has
\begin{align*}
\int_\Y\Big|\E[&f\!F_{y'}({X}_{k+1})|{ X}_k=\omega_k^i]-\E[f\!F_{y'}(\bar{ X}_{k+1})|\bar{ X}_k=\omega_k^i]\Big|dy' \\
& \leq L_{\!f}(\E[|{ X}_{k+1}-\bar{ X}_{k+1}|] +[\Phi]^2\E[| X_k-\bar{ X}_k|]).
\end{align*}
\end{lemma}

\noindent\textbf{Proof}
We split the integrand difference into 3 terms $A_i$, $1\leq i\leq 3$ with
\begin{align*}
A_1&=\E[f\!F_{y'}({\bar X}_{k+1})|{ \bar X}_k=\omega_k^i]-\E[f\!F_{y'}({ X}_{k+1})|\bar{ X}_k=\omega_k^i],\\
A_2&=\E[f\!F_{y'}({X}_{k+1})|{\bar X}_k=\omega_k^i]-\E[f\!F_{y'}({ X}_{k+1})|X_k=\bar{ X}_k=\omega_k^i],\\
A_3&=\E[f\!F_{y'}({X}_{k+1})|{ X}_k=\bar{ X}_k=\omega_k^i]-\E[f\!F_{y'}({ X}_{k+1})|{ X}_k=\omega_k^i].
\end{align*}
For the first term $A_1$, using  assumption \textbf{(H1)}  one obtains 
\begin{align*}
\int_\Y|A_1|dy'
&\leq L_{\!f}\E[|\bar{ X}_{k+1}-{ X}_{k+1}|\ |\bar{ X}_k=\omega_k^i].
\end{align*}
For the second term $A_2$, we use that $\bar{X}_k$ is the projection of ${ X}_k$ onto $\Omega_k$, so that one has the inclusion of $\sigma$-fields $\sigma(\bar{ X}_k)\subset\sigma({ X}_k)$. 
Note also that $\{ X_k=\omega_k^i\}=\{ X_k=\bar{ X}_k=\omega_k^i\}$ by definition of the projection. Thus, one has
\begin{align*}
A_2  &  =\E\Big[\E[f\!F_{y'}( X_{k+1})|{ X}_k]-\E[f\!F_{y'}( X_{k+1})|{ X}_k=\bar{ X}_k]\big| \bar{ X}_k=\omega_k^i\Big]\\
&=\E[{P_k\fF_{y'}(X_k)-P_k\fF_{y'}(\bar{ X}_k)}| \bar{ X}_k=\omega_k^i],
\end{align*}
and from similar derivations as in Lemma \ref{lem:Pflip} (recall that $ X_k$ and $\bar{ X}_k$ always have the same mode) we get
\begin{align*}
\int_\Y|A_2|dy'
&\leq [\Phi]^2L_{\!f}\E[|X_k-\bar{ X}_k|\ | \bar{ X}_k=\omega_k^i].
\end{align*}
Similarly, for the third term $A_3$, one has
\begin{align*}
A_3&=\left(\E[f\!F_{y'}( X_{k+1})|{ X}_k=\bar{ X}_k] -\E[f\!F_{y'}( X_{k+1})| X_k]\right)\1_{( X_k=\omega_k^i)}\\
&=\left(P_k\fF_{y'}( X_k)-P_k\fF_{y'}(\bar{ X}_k)\right)\1_{( X_k=\omega_k^i)}.
\end{align*}
Since on $\{ X_k=\omega_k^i\}$ one also has $\bar{ X}_k=\omega_k^i$ by projection,  one obtains $A_3=0$. 
Summing the three bounds and taking the expectation on both sides yields the result.  
\hfill$\Box$
%
\begin{lemma}\label{lem:Psi-Psibar}
Under assumption \textbf{(H1-3)}, for all $0\leq~k~\leq~N~-~1$ and $\thetabar\in\mathcal{P}(\Omega_k)$ one has
\begin{align*}
\int_\Y d(\bar\Psi_k(\thetabar,y'),\Psi_k(\thetabar,y'))dy'
  & \leq \fFbar \ \big(\E[|{ X}_{k+1}-\bar{ X}_{k+1}|]+[\Phi]^2\E[| X_k-\bar{ X}_k|]\big),
\end{align*}
for the distance between measures in $\mathcal{P}(\X)$.  
\end{lemma}
%
\noindent\textbf{Proof} 
Let $g\in BL_1(\X)$ { and $y\in\Y$}.  One has, with the notation of the proof of Lemma \ref{lem:Psi-dist}
\begin{align*}
\lefteqn{\int_\X g(\xi',y)d\bar\Psi_k(\thetabar,y')(\xi')-\int_\X g(\xi',y)d\Psi_k(\thetabar,y')(\xi')}\\
  &=\Big(\sum_{i=1}^{\ell_{k}}\!\E[f\!F_{y'}(X_{k+1})|X_k=\omega_k^i] \thetabar(\omega_k^i)\Big)^{-1} \\
&\qquad \times \Big(\sum_{i=1}^{\ell_{k}}\left(\E[f\!F_{y'}[g](\bar{X}_{k+1},y)|\bar{X}_k=\omega_k^i] \right. \left. -\E[f\!F_{y'}[g](X_{k+1},y)|X_k=\omega_k^i]\right)\!\thetabar(\omega_k^i) \Big)\\
&\quad +\sum_{i=1}^{\ell_{k}}\E[f\!F_{y'}[g](\bar {X}_{k+1},y)|\bar{X}_k=\omega_k^i]\thetabar(\omega_k^i)\\
&\qquad \times\left(\big(\sum_{i=1}^{\ell_{k}}\E[f\!F_{y'}(\bar{X}_{k+1})|\bar{X}_k=\omega_k^i]\thetabar(\omega_k^i)\big)^{-1} \right.  \left. -\big(\sum_{i=1}^{\ell_{k}}\E[f\!F_{y'}(X_{k+1})|X_k=\omega_k^i]\thetabar(\omega_k^i)\big)^{-1}\right)\\
&=A+B.
\end{align*}
By assumption \textbf{(H2)}, all denominators are bounded from below by $\underline{f}>0$.\\
$\bullet$ First term $A$: Following the same lines as the proof of Lemma \ref{lem:r-rbar integral}, one gets:
\begin{align*}
|A|\leq& 
\underline{f}^{-1}(L(y')+B(y'))\E[|\bar{X}_{k+1}-{X}_{k+1}|]+\underline{f}^{-1}[\Phi]^2(L(y')+B(y'))\E[|X_k-\bar{X}_k|].
\end{align*}
$\bullet$ Second term $B$: 
By similar computations as for term $A$, one has
\begin{align*}
B&\leq
 \frac{B(y')}{\underline{f}^2}L(y')\Big(\E[|\bar{X}_{k+1}-{X}_{k+1}|]\!+\![\Phi]^2\E[|X_k-\bar{X}_k|]\Big),
\end{align*}
hence the result.
 \hfill$\Box$
%
 \begin{lemma}\label{lem:Rk-Rkbar} 
Under assumption \textbf{(H1-3)}, for all $1\leq k\leq N$, $h~\in~ BLP(\X)$, and $(\thetabar,y)\in\mathcal{P}(\Omega_k)\times\Y$ one has
\begin{align*}
|\bar R'_kh(\thetabar,y)-R'_kh(\thetabar,y)|
\leq \Big( \|h\|_{\mathcal{P}}L_{\!f}+\overline{f}[h]_{\mathcal{P}}\fFbar\Big) 
(\E[|\bar{X}_{k+1}-{X}_{k+1}|]+[\Phi]^2\E[|X_k-\bar{X}_k|\ |]).
\end{align*}
Under assumption \textbf{(H2-3)}, for all $1\leq k\leq N$, $h\in BLP(\Omega_{k+1})$ not depending on $y$, and $(\thetahat,\yhat) \in\Gamma_k$  one has
\begin{align*}
{|\Rhat'_kh(\thetahat,\yhat) -\bar R'_kh(\thetahat,\yhat)|} 
&\leq \barkplip{h}\E\left[\dkp{\Thetahat_{k+1}-\Thetabar_{k+1}}\right] + \left(\barkpsup{h}+\fFtilde\overline{f} \barkplip{h}\right) \E\left[\dk{\Thetahat_{k}-\Thetabar_{k}}\right].
\end{align*}
\end{lemma}
 The proof of the first statement is a direct consequence of Lemmas \ref{lem:r-rbar integral} and \ref{lem:Psi-Psibar}. The proof of the second statement is obtained \textit{via} the same decomposition as that of Lemma \ref{lem:r-rbar integral} and is a consequence of Lemma \ref{Rk'-lip}.
%
\subsection{Proof of the main theorems}
\label{sec:proof-theo}
We are now able to prove our main theorems. Since both theorems can be proven using the same lines, we leave the proof of \ref{lem:vbark-vhatk} to the reader. While Theorem \ref{lem:v_n} requires the first statements of Lemma \ref{stability'} and Lemma \ref{lem:Rk-Rkbar}, the second theorem requires the second statements of those same lemmas. Recall also that none of the value functions depend on $y$.
%
\noindent\textbf{Proof of Theorem \ref{lem:v_n}}
Let $\theta$ be a distribution on $\X$ and $\thetabar$ on $\Omega_N$. Then, using Lemma \ref{lem:ga-BL}, for all $y,\ybar\in\mathbb{Y}$, one has 
  \begin{eqnarray*}
    |v'_N(\theta, y) - \vbar_N'(\thetabar,\ybar)|    &=& \left|\min_{a\in\mathcal{M}} C'(\theta,y,a) - \min_{a\in\mathcal{M}} C'(\thetabar,\ybar,a)\right|\\
    &\leq&\bar Bd(\theta,\thetabar).
  \end{eqnarray*}
Let $k< N$, $\theta$ a distribution on $\X$ and $\thetabar$ a distribution on $\Omega_k$. Then one has
\begin{align*}
|v'_k(\theta, y) - \vbar_k'(\thetabar,\ybar)|
  & \leq |v'_k(\theta, y)-v'_k(\thetabar, \ybar)|+|v'_k(\thetabar, \ybar)-\vbar_k'(\thetabar,\ybar)|.
\end{align*}
The first term is readily bounded using Lemma \ref{stability'}. 
For the second term, one has 
  \begin{align*}
|v'_k(\bar\theta, \ybar) - \vbar_k'(\thetabar,\ybar)|
    & \leq   \left|\Rbar'_k\vbar'_{k+1}(\thetabar,\ybar)-\Rbar'_kv'_{k+1}(\thetabar,\ybar)\right|  +   \left|\Rbar'_kv'_{k+1}(\thetabar,\ybar)-R'_kv'_{k+1}(\thetabar,\ybar)\right|. 
  \end{align*}
We study the two terms above separately.  
For the first term, as $\Rbar'_k$ is a Markov kernel and according to Lemma \ref{stability'}, $v_{k+1}' \in BLP(\X)$, one clearly has $\left|\Rbar'_k(\vbar'_{k+1} -v'_{k+1})(\thetabar,\ybar)\right|\leq \|\vbar'_{k+1} -v'_{k+1}\|_{k+1}$ (note that the $\sup$ is taken on $\Omega_{k+1}$).
For the second term, we use Lemma \ref{lem:Rk-Rkbar} to obtain
  \begin{align*}
\left|\Rbar'_kv'_{k+1}(\thetabar,\ybar)-R'_kv'_{k+1}(\thetabar,\ybar)\right|
    &\leq  \big( L_{\!f}\|v'_{k+1}\|_{\mathcal{P}}+\fFbar[v'_{k+1}]_{\mathcal{P}\overline{f}}\big) 
(\E[|\bar{X}_{k+1}-{X}_{k+1}|]+[\Phi]^2\E[|X_k-\bar{X}_k|]).
  \end{align*}
Combining the previous bounds leads to
the result.
\hfill$\Box$
%
\section{Simulation study}
\label{Simulations}
%
In this section we consider several simulation studies to assess the performance of our candidate strategy on the models of Section \ref{sec:examples}. In the simplest scenario, we compare it to other state of the art methods. Recall that the best (possibly unfeasible) performance is $0$.

The performance of the approaches is evaluated through the average cost of the strategies over $1000$ Monte Carlo simulations. We consider different combinations of values for the costs of stopping too early, stopping too late, or choosing the wrong mode, namely
\begin{align*}
  \alpha & \in \{3,4,5,6\} \\
  \beta& \in \{0.5,1,1.5,2\} \\
  \gamma_{ma}= \gamma  & \in \{0.5,1,1.5,2\}.
\end{align*}
In all simulations we consider a centered Gaussian noise $\varepsilon_n$ for the observations with three possible variance values : $0.1, 0.5$ and $1$. The finite real-time horizon is set to $\mathcal{T}=6$, so that $\mathbb{P}(T>\mathcal{T}) \simeq 1.5e-08$ and the optimization horizon $N$ is~$\left \lfloor{\frac{\mathcal{T}}{\delta}}\right \rfloor$, with $\delta$ set to $1/6$.
For the quantization grids, we fixed $\ell_k$ the number of points in $\Omega_k$ to $21$ for all observation times, and considered $4$ values for $N_k$ the number of points in $\Gamma_k$ (once again fixed for all $k$): $30,50,75$ and $100$. The grids were calibrated using the CLVQ algorithm presented for instance in \cite{pages04b}. 
%
\subsection{Models \ref{Ex:expo}}
\label{sec:comparison}
%
In this section we first compare our method with two other state of the art approaches, namely a moving average approach and  Kalman filtering,  using  models presented in Example \ref{Ex:expo}.

\paragraph*{Moving average (MA)} is a model-free online strategy that has been widely used to detect change-points in time series. Specifying a window size $k$ and a threshold $s$, the method consists, at time $n$, in computing the mean $T_{n}^k$ of the signal over observations $\{y_{n-k+1}, y_{n-k+2}, \dots, y_{n}\}$ and declaring a change point if  $T_{n}^k>s$  or, depending on the challenge at hand,  $T_{n}^k<s$. \\
\underline{Parameter tuning:}
we consider $4$ possible window sizes $k$ taking values in $\{2,3,4,5\}$ and $3$ possible thresholds $s$: $\{1.25,1.5,2\}$. \\ 
\underline{Mode-choice rule:}
we propose a heuristic method to select the mode after the detection is performed. Assuming the detection occurs at observation time $\tau$, for each possible $v_i$, we compute the minimum over $t \in \{2, \dots \tau-1\}$ of the inverse log-likelihood $\ell_i$ of the data under the (wrong) assumption that the $Y_n$ are independent, following an $\mathcal{N}(1,\sigma^2)$ up to $t$ and an $\mathcal{N}(e^{v_i(n-t)\delta},\sigma^2)$ for $n$ between $t+1$ and $\tau$. We  then select the mode $i$ which minimizes the $\ell_i$.
\paragraph*{Kalman filtering (KF)} is a model-based online strategy that estimates the states of a system in the particular case of linear dynamical systems. In our simulation scenario, our model can be approximated to fit  a Kalman model in the following way :
\begin{equation*}
  \left\{
  \begin{array}{rl}
    {X}_n&=A_{M_n} X_{n-1}, \\
    {Y}_n&=X_n + \varepsilon_n, \   \varepsilon_n \sim \mathcal{N}(0,\sigma^2); \\
  \end{array}
  \right.
\end{equation*}
%
with $M_n$ a Markov chain on $\{0,1,\dots,d\}$ that represents the mode of ${X}_n$ at time $t_n$ and $A_i=e^{v_i\delta}$ with the convention $v_0=0$. Of note, this approximation would be exact if $T$ occurred on an observation date. \\
\underline{Transition matrix:}
applying KF requires the knowledge of the mode transition matrix and can easily be extended to inhomogeneous Markov chains. In our simulation scenario, this transition matrix is computed analytically. \\
\underline{Threshold calibration:}
this approach also requires the choice of a threshold $s$ for which if $\mathbb{P}(M_i|Y_1,\dots Y_\tau)>s$, the decision to stop the process and select mode $i$ is taken. In this simulation study we consider $4$ possible values for $s$ : $0.5, 0.75$, $0.9$ and a calibrated threshold based on the costs and defined as follows.
At observation time $k$, the integrated cost of taking decision $a=\Delta_j$ is  $C^k_j= C^{k-1}_0 + \sum_{i=0}^d u(\Delta_j|m=i)\mathbb{P}(M_i|Y_1,\dots Y_k)$, and one should take the decision with minimal cost (with convention $C^0_0=0$). Hence in practice, the following decision rule can be followed: \\
$\bullet$ if for all $1\leq j \leq d$, $\beta \delta \sum_{i=1}^d \mathbb{P}(M_i|Y_1,\dots Y_k) \leq \alpha \mathbb{P}(M_0|Y_1,\dots Y_k) + \gamma \sum_{i\neq j} \mathbb{P}(M_i|Y_1,\dots Y_k)$ then continue the observation process until the next observation. \\
$\bullet$ else $\tau=k$ and take decision $j$ that minimizes $ \alpha~\mathbb{P}(M_0|Y_1,\dots Y_k)~ +~ \gamma~ \sum_{i\neq j} \mathbb{P}(M_i|Y_1,\dots Y_k)$.

\paragraph*{Comparison for Model \ref{Ex:expo}.a} Here we selected $d=3$ and $v_i \in \{0.1, 0.5, 1\}$. Table \ref{tab:expo} presents the results for $\beta=1$ (cost of stopping too late), $\gamma=1.5$ (cost of selecting the wrong mode) and the four values of $\alpha$ (cost of an early detection). For MA, only the results for a threshold value of $2$ are presented as they are always better. 
\begin{table}[tp]
  \centering
      \scalebox{.9}{
        \begin{tabular}{|r|r|cccc|cccc|cccc|}
        \hline
          &   & \multicolumn{4}{c|}{MA, $s=2$}   & \multicolumn{4}{c|}{KF} & \multicolumn{4}{c|}{New approach}\\
          \hline
          $\alpha$ &$\sigma^2$      & \multicolumn{4}{c|}{window}    & \multicolumn{4}{c|}{threshold}  & \multicolumn{4}{c|}{$N_k$} \\
          \hline
          &&  2 & 3 & 4 & 5  & 0.5 & 0.75 & 0.9  & cal & 30 & 50 & 75 & 100 \\ 
          \hline
          &0.1 &  0.41 & \textcolor{blue}{\textbf{0.40}} & \textcolor{blue}{\textbf{0.40}} & 0.41 & 2.30 & 0.59 & 0.42 & \textbf{0.41} & 0.70 & \textbf{0.69} & 0.70 & \textbf{0.69} \\ 
          3&0.5 &  0.92 & 0.81 & 0.76 & \textbf{0.72} & 1.45 & 0.56 & 0.49 & \textcolor{blue}{\textbf{0.48}} & \textbf{0.79} & 0.84 & 0.83 & 0.82 \\ 
          &1 &  1.70 & 1.45 & 1.27 & \textbf{1.18} & 1.22 & \textcolor{blue}{\textbf{0.58}} & 0.65 & 0.66 & \textbf{0.99} & 1.00 & 1.00 & 1.00 \\
          \hline
          \hline
          &0.1 &  0.41 & \textcolor{blue}{\textbf{0.40}} & \textcolor{blue}{\textbf{0.40}} & 0.41 & 3.00 & 0.66 & 0.42 & \textbf{0.41} & 0.68 & 0.69 & \textbf{0.67} & 0.68 \\ 
          4&0.5 &  0.93 & 0.81 & 0.76 & \textbf{0.72} & 1.78 & 0.59 & 0.49 & \textcolor{blue}{\textbf{0.48}} & 0.75 & 0.75 & \textbf{0.74} & \textbf{0.74} \\ 
          &1 &  2.02 & 1.61 & 1.35 & \textbf{1.23} & 1.42 & \textcolor{blue}{\textbf{0.59}} & 0.65 & 0.66 & 0.96 & 0.90 & \textbf{0.89} & 0.91 \\
          \hline
          \hline
          &0.1 &  0.41 & \textcolor{blue}{\textbf{0.40}} & \textcolor{blue}{\textbf{0.40}} & 0.41 & 3.70 & 0.74 & 0.42 & \textbf{0.41} & 0.67 & \textbf{0.64} & 0.66 & 0.66 \\ 
          5&0.5 &  0.94 & 0.81 & 0.77 & \textbf{0.73} & 2.11 & 0.62 & 0.49 & \textcolor{blue}{\textbf{0.48}} & \textbf{0.70} & 0.71 & 0.72 & 0.73 \\ 
          &1 & 2.34 & 1.78 & 1.44 & \textbf{1.29} & 1.62 & \textcolor{blue}{\textbf{0.60}} & 0.65 & 0.66 & 0.96 & 0.90 & \textbf{0.89} & 0.90 \\
          \hline
          \hline
          &0.1 &  0.41 & \textcolor{blue}{\textbf{0.40}} & \textcolor{blue}{\textbf{0.40}} & 0.41 & 4.41 & 0.82 & 0.42 & \textbf{0.41} & 0.66 & \textbf{0.64} & 0.67 & 0.65 \\ 
          6&0.5 &  0.95 & 0.82 & 0.77 & \textbf{0.73} & 2.44 & 0.65 & 0.49 & \textcolor{blue}{\textbf{0.48}} & \textbf{0.67} & 0.69 & 0.70 & 0.70 \\ 
          &1 &  2.67 & 1.94 & 1.52 & \textbf{1.35} & 1.82 & \textcolor{blue}{\textbf{0.61}} & 0.65 & 0.66 & 0.93 & 0.88 & \textbf{0.86} & 0.87 \\
          \hline
        \end{tabular}
      }
      \caption{{ Average cost of the strategies for MA, KF and the approach proposed in this paper for different parameter values.} Best scores  (\textit{i.e.} minimal strategy cost) by method are indicated in bold, while the overall best score is indicated in blue.} \label{tab:expo}
    \end{table} 
Though MA performs very well for small variance values, it gets out-beaten by KF for intermediate variance values, and by both KF and our approach for high variance values, with optimal performance for large window and threshold.  Moreover, the impact of increasing the noise value is tremendous for MA, while both other approaches adapt more easily. KF out-beats the others as soon as the variance is intermediate (with a performance almost as good as MA for the low variance scenario), but the optimal threshold depends highly on the noise level. Finally, our approach is better than MA for medium variance values, but is beaten by the KF. Interestingly, increasing the number of points in the quantization grid for the filter does not improve the performance of the method. This is expected with small $N_k$ increments since the quantization convergence rate is in $\mathcal{O}(N_k^{-1/\ell_k}) =\mathcal{O}(N_k^{-1/21})$.
Finally, none of the methods are significantly affected by the increase in penalty value for early detection. This comes from the fact that, on this example, all methods tend to detect the mode jump with delay, the earliest detection being performed by KF with probability threshold $0.5$.

As seen above, KF out-beats the other methods. This is expected since our simulation scenario is very close to a linear system for which KF is known to be optimal. However, when the observation frequency decreases, for instance for $\delta=1/2$, the Kalman model approximation of our simulation deteriorates. This leads to degradation of the KF  performances, and to our approach becoming optimal. 
We present two other types of deviation from this model, on which we compare our PDMP approach to Kalman filtering.
\paragraph*{Comparison for Model \ref{Ex:expo}.a with inverse link function} We use Model \ref{Ex:expo}.a again, this time replacing the link function $F(x)=x$ by  $F(x)=1/x$. 
To apply KF, we consider a linearized version of the observations, which actually corresponds to the same Kalman model as above. 
Table \ref{tab:expoF} presents the results for $\sigma^2=0.5$, a filter quantization grid size of $50$ points, $\alpha=4$, $\beta=1$ and the four values of $\gamma$. Here our approach has better results than KF in terms of strategy cost. Moreover it adapts more easily to the increase in penalty value for selecting the wrong mode, suggesting that the difference in performance is mostly driven by the mode decision.
\begin{center}
  \begin{table}[tp]
    \centering
    \scalebox{.9}{
    \begin{tabular}{c|cccc|c|}
      &\multicolumn{4}{c|}{Kalman} & New approach \\
      \hline
      $\gamma$ & 0.5 & 0.75 & 0.9 & cal & $N_k=50$  \\ 
      \hline
      0.50 & 1.55 & \textbf{0.75} & 0.98 & 0.86 & \textcolor{blue}{\textbf{0.69}} \\ 
      1.00 & 1.78 & \textbf{1.07} & 1.34 & 1.14 & \textcolor{blue}{\textbf{0.93}} \\ 
      1.50 & 2.02 & \textbf{1.40} & 1.69 & 1.65 & \textcolor{blue}{\textbf{1.11}} \\ 
      2.00 & 2.26 & \textbf{1.73} & 2.04 & 2.18 & \textcolor{blue}{\textbf{1.32}} \\ 
      \hline
    \end{tabular}}
    \caption{{ Average cost of the strategies for  KF and our approach on a non-trivial link function simulation study}}\label{tab:expoF}
  \end{table} 
\end{center}
\paragraph*{Comparison for Model \ref{Ex:expo}.b}
To apply the KF approach to Model \ref{Ex:expo}.b, we consider a first order approximation of the exponential model, so that $A_1=e^{v_1\delta}$ and $A_2=v_2\delta$. Values used are $v_1=0.5$ for the exponential flow, and $v_2=2$ for the linear flow.
Table \ref{tab:expoDroite} presents the results for $\sigma^2=0.5$,  $\alpha=4$, $\gamma=1.5$ and the four values of $\beta$. Here despite the fact that one of the possible modes corresponds to a linear model, KF fails to find an optimal strategy: the algorithm leads to early detection in 53-69 \% of cases for the calibrated decision rule, and as much as 90\% of cases for the 0.5 threshold. In comparison, our approach leads to early detection in 0\% (for the lowest $\beta$ value) to 12 \% (for the highest $\beta$ value) of cases, explaining the large difference in strategy cost.
\begin{center}
  \begin{table}[tp]
    \centering
    \scalebox{.9}{
    \begin{tabular}{c|cccc|ccc|}
      &\multicolumn{4}{c|}{Kalman} & \multicolumn{3}{c|}{New approach} \\
      \hline
       $\beta$ & 0.5 & 0.75 & 0.9 & cal & 20 & 100 & 200 \\ 
      \hline
      0.50 & 3.38 & 3.01 & 2.56 & \textbf{2.11} & 0.38 & \textcolor{blue}{\textbf{0.30}} & \textcolor{blue}{\textbf{0.30}} \\ 
      1.00 & 3.38 & 3.01 & 2.57 & \textbf{2.33} & 0.52 & 0.50 & \textcolor{blue}{\textbf{0.49}} \\ 
      1.50 & 3.38 & 3.02 & 2.58 & \textbf{2.51} & 0.70 & 0.71 & \textcolor{blue}{\textbf{0.65}} \\ 
      2.00 & 3.38 & 3.02 & \textbf{2.58} & 2.66 & 0.86 & 0.89 & \textcolor{blue}{\textbf{0.85}} \\ 
      \hline
    \end{tabular} }
    \caption{{ Average cost of the strategies for  Kalman and PDMP on a non-linear model}}\label{tab:expoDroite}
  \end{table} 
\end{center}
%
\subsection{Models \ref{Ex:sin1}}
We turn to Example \ref{Ex:sin1} for which, to our knowledge, no other algorithm is adapted to change-point detection. For Model \ref{Ex:sin1}.a, we set $d=1$ and $v_0=3, v_1=5.$ For Model \ref{Ex:sin1}.a we set $d=2$ and $v_0=3, v_i \in \{0.5,1.5\}$.  We considered a fixed number of points per quantization grid (respectively $10$ and $15$ for $\Omega_k$ in the a and b scenarios, and $50$ for $\Gamma_k$ in both cases).
\begin{table}[tp]
  \centering
      \scalebox{.9}{
        \begin{tabular}{c|c|cccc|cccc|}
          &   & \multicolumn{4}{c|}{$\alpha=3$, $\beta=2$}   & \multicolumn{4}{c|}{$\alpha=6$, $\beta=0.5$} \\
          \hline
          $\delta$  & $\sigma^2$ &  $\Delta T$ & sd & Nb Obs & Nb early  & $\Delta T$ & sd & Nb Obs & Nb early \\
          \hline
          &0.1 & 0.35 &0.06  &  4 &   0& 0.39& 0.06 &   4 &   0 \\
          1/10 &  0.5 & 0.48 & 0.17 &   5 &  16& 0.59 & 0.19  &  6 &  10 \\
          &1.0 & 0.44 & 0.44 &   4 & 103 & 0.67 & 0.46  &  7 &  60 \\
          \hline
          \hline
          & 0.1& 0.44& 0.09 &   3 &   0& 0.54& 0.12 &   3  &  0 \\
          1/6 &  0.5& 0.59& 0.23 &   4  & 20& 0.71& 0.26 &   4 &  12 \\
          & 1.0& 0.48& 0.51 &   3 & 124&0.91& 0.54 &   5 &  38 \\
          \hline
          \hline
          &0.1& 0.71& 0.48  &  3 &  10& 1.06 &0.95   & 4   & 0 \\
          1/4 & 0.5 & 0.67 & 0.54 &   3 &  81&1.27& 0.78   & 5   & 3 \\
          &1.0 & 0.57 & 0.62 &   2 & 152& 1.28& 0.73 &   5  & 22 \\
          \hline
        \end{tabular}
      }
      \caption{{Time to jump detection  for different $\delta$ values.} Average real-time and observation-time to detection of our approach on $1000$ simulations from Model~\ref{Ex:sin1}.a) } \label{tab:timesin}
    \end{table} 
Table \ref{tab:timesin} shows the influence of the observation time-steps $\delta$ and of the variance $\sigma^2$ in terms of time to jump-detection for Model~\ref{Ex:sin1}.a (results for Model~\ref{Ex:sin1}.b are very similar). While increasing $\delta$ significantly increases the amount of time required to detect the jump, it slightly decreases the number of observations needed after the jump for its detection. The number of early detection is more influenced by the cost parameters than by the time steps, with a strong tendency to increase with the observation noise. This in fact leads to strategy costs very close to zero for small time-steps and variance values, with a tendency to increase with $\delta, \sigma^2$ and to decrease with $\alpha$.
%
\section{Conclusion}
\label{sec:CCL}
%
We have proposed a numerically feasible numerical scheme to approximate the value function of a change-point detection problem for a simple class of PDMPs.
We obtain error bounds for this approximation explicitly depending on the parameters of the problem. We have also proposed a feasible stopping strategy that performs well compared to state of the art methods when such methods are applicable, despite the very slow convergence rate.
We believe that this is a promising start for the more general study of impulse control problems for general PDMPs when there is no observations of the jump times. The easiest extension is certainly going from scalar-valued PDMPs to multivariate ones. We believe that our proofs would hold in this context.
Allowing more than one jump should be more challenging as one would not be able to write the explicit form of the kernels $P_n$. However, the underlying POMDP framework is suitable for several interventions so that it should be possible to extend our results in this direction, although probably technically involved. 
Finally, the important open question concerns the optimality of our candidate strategy. It cannot be directly linked to our various operators, but we are hopeful that further work will enable us to prove theoretically that it is close to optimality.
\appendix
\section{Proofs}
\label{sec:appx}
%
\begin{lemma}\label{add'}
For any function $g\in BL(E)$ set $g':(\theta,y)\mapsto \int g(\xi,y)d\theta(\xi)$ for $(\theta,y)\in \mathcal{P}(E)\times \Y$. Then $g'$ is in $BLP(E)$ with
\begin{align*}
  \|g'\|_{\mathcal{P}}   \leq\|g\|,\quad  {[g']}_{\mathcal{P}}  &\leq  \|g\|+[g].
\end{align*}
\end{lemma}
%
\noindent\textbf{Proof}
Clearly, one has
$|g'(\theta,y)|\leq \int_\X\|g\|d\theta(x)=\|g\|$ and 
\begin{align*}
|g'(\theta,y)-g'(\thetabar,y)|
&\leq \Big|\int_\X g(m,x,y)d\theta(m,x)-\int_\X g(m,x,y)d\thetabar(m,x)\Big|\\
&\leq \Big|\int_\X g(m,x,y)d\theta(m,x)-\int_\X g(m,x,y)d\thetabar(m,x)\Big|\\
&\leq ([g]+\|g\|)d_E(\theta,\thetabar),
\end{align*}
hence the result.
  \hfill$\Box$\\

\noindent\textbf{Proof of Lemma \ref{lem:ga-BL}}
 One has 
 $$C_a(m)=\alpha\1_{m=0}+\gamma_{m,a} \1_{m>0,\neq a},$$
 thus it is easy to see that $\|C_a\|= \max\{\alpha, \{\gamma_{ia}\}_{1\leq i \leq d}\}$ and as $C_a$ does not depend on $x$ nor $y$, $[C_a]=0$.
 Similarly, $c(m,x,y)=\beta\delta\1_{m\neq 0}$ thus $\|c\|=\beta\delta$ and $[c]=0$.
The other statements follow from an application of Lemma \ref{add'}.
  \hfill$\Box$
%
\begin{lemma}\label{fFg}
For any function $g\in BL(\X)$ and $y'\in\Y$, the function $f\!F_{y'}[g]:(m,x,y)\mapsto g(m,x,y)f(y'-F(x))$ is in $BL(\X)$ with
\begin{align*}
  \|f\!F_{y'}[g]\|  \leq \|g\|B(y'),
  \quad  {[f\!F_{y'}[g]]} \leq \|g\|L(y')+ [g]B(y').
\end{align*}
\end{lemma}
%
\noindent\textbf{Proof}
By definition, one has $\|f\!F_{y'}[g]\|\leq \|g\|\sup_{x\in\K}f(y'-F(x))$ and 
\begin{align*}
\lefteqn{|g(m,x,y)f(y'-F(x))-g(m,x',y)f(y'-F(x'))|}\\
&\leq |g(m,x,y)f(y'-F(x))-g(m,x,y)f(y'-F(x'))|+|g(m,x,y)f(y'-F(x'))-g(m,x',y)f(y'-F(x'))|\\
&\leq \|g\|L(y')|x-x'|+[g]\sup_{x\in\K}f(y'-F(x))|x-x'|,
\end{align*}
hence the result.
  \hfill$\Box$\\

\noindent\textbf{Proof of Lemma \ref{lem:Psi-dist}, second part}
The proof follows the same lines as that of Lemma \ref{lem:Psi-dist}. By definition, one has
  \begin{align*}
    \lefteqn{\dkp{\Psibar_{k}(\theta,y')-\Psibar_{k}(\thetabar,y')}}\\
    &\leq    \sum_{j=1}^{\ell_{k+1}}\left|\frac{\sum_{i=1}^{\ell_k}\bar p_{k,i,j}f(y'-F(\omega_{k+1}^j))\big(\theta(\omega_k^i)-\thetabar(\omega_k^i)\big)}{\sum_{j'=1}^{\ell_{k+1}}\sum_{i=1}^{\ell_k}\bar p_{k,i,j}f(y'-F(\omega_{k+1}^{j'}))\theta(\omega_k^i)}\right| \\
    & \quad + \sum_{j=1}^{\ell_{k+1}}\sum_{i=1}^{\ell_k}\bar p_{k,i,j}f(y'-F(\omega_{k+1}^j))\thetabar(\omega_k^i) \\
    & \qquad \left|\frac{1}{\sum_{j'=1}^{\ell_{k+1}}\sum_{i=1}^{\ell_k}\bar p_{k,i,j}f(y'-F(\omega_{k+1}^{j'}))\theta(\omega_k^i)} \right. 
 \left.-\frac{1}{\sum_{j'=1}^{\ell_{k+1}}\sum_{i=1}^{\ell_k}\bar p_{k,i,j}f(y'-F(\omega_{k+1}^{j'}))\thetabar(\omega_k^i)}\right|\\
    &= A +B.
  \end{align*}
  As in the first part of the proof of Lemma \ref{lem:Psi-dist}, from assumption {\bf (H2)} we have all denominators bounded from below by $\underline{f}$.
  For the first term $A$, we therefore have
  \begin{align*}
 A   & \leq \frac{1}{\underline{f}}\sum_{j=1}^{\ell_{k+1}}\sum_{i=1}^{\ell_k} \bar p_{k,i,j}f(y'-F(\omega_{k+1}^j))\left|\theta(\omega_k^i)-\thetabar(\omega_k^i)\right| \\
    & \leq 
    \frac{B(y')}{\underline{f}} \dk{\theta-\thetabar}.
  \end{align*}
  For the second term $B$, we have
  \begin{align*}
    B   & \leq \frac{\nsup{\fF}}{\underline{f}^2}\Big|{ \sum_{j'=1}^{\ell_{k+1}} \sum_{i=1}^{\ell_k} \bar p_{k,i,j}f(y'-F(\omega_{k+1}^{j'})) \big(\theta(\omega_k^i)- \thetabar(\omega_k^i)\big)}\Big| \\
    & \leq  \frac{B(y')^2}{\underline{f}^2} \dk{\theta-\thetabar},
  \end{align*}
  hence the result.
\hfill$\Box$\\

\noindent\textbf{Proof of Lemma \ref{Rk'-lip}, second part}
 The first statement is obvious as $\Rbar'_{k}$ is a Markov kernel.
 Now let $(\theta,\thetabar) \in \mathcal P(\Omega_{k})^2$ and $(y,\ybar) \in \Y^2$.
 \begin{align*}
   \lefteqn{\left| \Rbar'_{k} g(\theta,y)- \Rbar'_{k} g(\thetabar,\ybar) \right|}\\
&\leq \int_{\Y}\left|g(\Psibar_{k}(\theta,y'),y')-g(\Psibar_{k}(\thetabar,y'),y')\right| 
\sum_{\substack{1\leq i\leq \ell_{k}\\1\leq j\leq \ell_{k+1}}}\!\!\!\!\!\bar p_{n,i,j}f(y'-F(\omega_{k+1}^j))\theta(\omega_k^i)dy'\\
&+\int_{\Y}\!g(\Psibar_{k}(\thetabar,y'),y')\! \!\!\!\!\!\!\sum_{\substack{1\leq i\leq \ell_{k}\\1\leq j\leq \ell_{k+1}}}\!\!\!\!\!\!\bar p_{k,i,j}f(y'-F(\omega_{k+1}^j))|\theta(\omega_k^i)-\thetabar(\omega_k^i)| dy'\\
&\leq ( [g]_{k+1}\overline{f} \fFtilde+\|g\|_{k+1})\|\theta -\thetabar\|_{k,1},
 \end{align*}
 using lemma \ref{lem:Psi-dist}.
\hfill$\Box$\\

\noindent\textbf{Proof of Lemma \ref{stability'}, second part}
 We proceed by backward induction on k.\\
 $\bullet$ For $k=N$, let $(\theta,y)$ and $(\thetabar,y) \in \mathcal P(\Omega_{N})\times \Y$.
 From Lemma \ref{lem:ga-BL}, on the one hand, one has
 \begin{align*}
   \left| \vbar'_{N} (\theta,y) \right| &= \big|\min_{a\in\mathbb A} C'(\theta,y,a) \Big| \leq  \delta\bar\beta\vee\bar G\leq \bar B, 
 \end{align*}
 and on the other hand,
 \begin{align*}
  \left| \vbar'_{N} (\theta,y)-\vbar'_{N} (\thetabar,y) \right| 
   &= \big|\min_{a\in\mathbb A} C'(\theta,y,a)-\min_{a\in\mathbb A} C'(\thetabar,y,a) \big| \leq  \bar B \left\|\theta-\thetabar\right\|_{{N,1}}.
 \end{align*}
 $\bullet$ We now assume that the statement holds for a given $k+1$ in $\{1,\dots, N\}$.
 Let $\theta \in \mathcal P(\Omega_{k})$ and $y \in \Y$. According to Lemma \ref{Rk'-lip}, $\vbar'_{k+1}$ is in $BLP(\Omega_{k+1})$, thus
 \begin{eqnarray*}
   \barksup{\vbar'_k} &\leq&\max_{a\in\mathcal{M}_+} \barksup{C_a}\vee\big( \barksup{c'}+ \barkpsup{\vbar'_{k+1}}\big)\\
   &\leq&
    = \bar B  + \barkpsup{\vbar'_{k+1}} \\
   &\leq & (N-k+1)\bar B.
 \end{eqnarray*}
 Regarding the Lipschitz property, Lemma \ref{Rk'-lip} yields 
 \begin{eqnarray*}
 \big|\vbar'_k(\theta,y)-\vbar'_k(\thetabar,y)\big|
   &\leq&\left(\bar B  +\barkpsup{\vbar'_{k+1}}+\barkplip{v'_{k+1}}\fFtilde\overline{f} \right) \ \dk{\theta-\thetabar},
  \end{eqnarray*} 
 hence the result.
\hfill$\Box$\\

\noindent\textbf{Detailed proof of Lemma \ref{lem:Psi-Psibar}}
To find an upper bound for $A$,
we follow the same lines as in the proof of Lemma \ref{lem:r-rbar integral}:
 \begin{align*}
 \lefteqn{  \E[f\!F_{y'}[g](\bar{X}_{k+1},y)|\bar{X}_k=\omega_k^i]
 -\E[f\!F_{y'}[g](X_{k+1},y)|X_k=\omega_k^i]}\\
 & =   \E[f\!F_{y'}[g](\bar{X}_{k+1},y)|\bar{X}_k=\omega_k^i]
 -  \E[f\!F_{y'}[g]({X}_{k+1},y)|\bar{X}_k=\omega_k^i] \\
 &\quad +  \E[f\!F_{y'}[g]({X}_{k+1},y)|\bar{X}_k=\omega_k^i]
 -  \E[f\!F_{y'}[g]({X}_{k+1},y)|{X}_k=\bar{X}_k=\omega_k^i] \\
 &\quad +   \E[f\!F_{y'}[g]({X}_{k+1},y)|{X}_k=\bar{X}_k=\omega_k^i]
 -  \E[f\!F_{y'}[g]({X}_{k+1},y)|{X}_k=\omega_k^i] \\
 &=A_1+A_2+A_3.
 \end{align*}
 For the first term $A_1$, assumption \textbf{(H3)} and Lemma \ref{fFg} yield
 \begin{align*}
 |A_1|
 &\leq (\|g\|L(y')+ [g]B(y'))\E[|\bar{X}_{k+1}-{X}_{k+1}|\ |\bar{X}_k=\omega_k^i].
 \end{align*}
 For the second term $A_2$, similarly to the proof of Lemma \ref{lem:r-rbar integral}, we use the properties of the projection and Lemma~\ref{lem:Pflip} to obtain 
 \begin{align*}
 |A_2|
 &\leq [\Phi]^2(\|g\|L(y')+ [g]B(y'))\E[|X_k-\bar{X}_k|\ | \bar{X}_k=\omega_k^i].
 \end{align*}
 Similarly, the third term $A_3$ equals $0$. 
 Summing the three terms, one gets
 \begin{align*}
 {|A|}
 &\leq \frac{1}{\underline{f}}\sum_{i=1}^{\ell_{k}}\thetabar(\omega_k^i)(|A_1|+|A_2|+|A_3|)\\
 &\leq \frac{(\|g\|L(y')+ [g]B(y'))}{\underline{f}}\sum_{i=1}^{\ell_{k}}\thetabar(\omega_k^i)(\E[|\bar{X}_{k+1}-{X}_{k+1}|\ |\bar{X}_k=\omega_k^i]+ [\Phi]^2\E[|X_k-\bar{X}_k|\ | \bar{X}_k=\omega_k^i]).
 \end{align*}
 Taking the expectation on both sides yields
\hfill$\Box$\\

\noindent\textbf{Proof of Lemma \ref{lem:Rk-Rkbar}}
\textbf{First statement}
 One has 
\begin{align*}
 \lefteqn{\left|\bar R'_kh(\thetabar,y)-R'_kh(\thetabar,y)\right|} \\
 =&\Big|\int_\Y h(\bar\Psi_k(\thetabar,y'),y')\sum_{i=1}^{\ell_k}\E[f\!F_{y'}(\bar X_{k+1})|\bar X_k=\omega_k^i]\thetabar(\omega_k^i)dy'
-\int_\Y h(\Psi_k(\thetabar,y'),y')\sum_{i=1}^{\ell_k}\E[f\!F_{y'}(X_{k+1})| X_k=\omega_k^i]\thetabar(\omega_k^i)dy'\Big|\\
\leq&\int_\Y |h(\bar\Psi_k(\thetabar,y'),y')|\sum_{i=1}^{\ell_k}\Big|\E[f\!F_{y'}(\bar X_{k+1})|\bar X_k=\omega_k^i]
-\E[f\!F_{y'}(X_{k+1})| X_k=\omega_k^i]\Big|\thetabar(\omega_k^i)dy'\\
&+\int_\Y \Big|h(\bar\Psi_k(\thetabar,y'),y')-h(\Psi_k(\thetabar,y'),y')\Big|
\sum_{i=1}^{\ell_k}\E[f\!F_{y'}(X_{k+1})| X_k=\omega_k^i]\thetabar(\omega_k^i)dy'
 \end{align*} 
and we conclude using lemmas \ref{lem:r-rbar integral} and \ref{lem:Psi-Psibar}.\\
\textbf{Second statement}
 We follow again the same lines as in the proof of Lemma \ref{lem:r-rbar integral}. One has 
 \begin{align*}
   \lefteqn{\Rhat'_kh(\thetahat,\yhat) - \bar R'_kh(\thetahat,\yhat)} \\
   & =\E\left[h(\Thetahat_{k+1},\Yhat_{k+1})|(\Thetahat_{k},\Yhat_{k})=\gamma_k^p\right]
   -\E\left[h(\Thetabar_{k+1},\Ybar_{k+1})|(\Thetahat_{k},\Yhat_{k})=\gamma_k^p\right]  \\
   & +\E\left[h(\Thetabar_{k+1},\Ybar_{k+1})|(\Thetahat_{k},\Yhat_{k})=\gamma_k^p\right]
   -\E\left[h(\Thetabar_{k+1},\Ybar_{k+1})|(\Thetabar_{k},\Ybar_{k})=(\Thetahat_{k},\Yhat_{k})=\gamma_k^p\right]  \\
   &  +\E\left[h(\Thetabar_{k+1},\Ybar_{k+1})|(\Thetabar_{k},\Ybar_{k})=(\Thetahat_{k},\Yhat_{k})=\gamma_k^p\right]  
   -\E\left[h(\Thetabar_{k+1},\Ybar_{k+1})|(\Thetabar_{k},\Ybar_{k})=\gamma_k^p\right]  \\
   & = A_1+ A_2 + A_3
 \end{align*} 
 $\bullet$  For the first term $A_1$, since $h \in BLP(\Omega_{k+1})$ and does not depend on $y$, one has
 \begin{align*}
   |A_1|
   & \leq \barkplip{h} \E\left[\dkp{\Thetahat_{k+1}-\Thetabar_{k+1}}|\big|(\Thetahat_{k},\Yhat_{k})=\gamma_k^p\right].
 \end{align*} 
 %
 $\bullet$  For the second term $A_2$, we use that $(\Thetahat_k,\Yhat_k)$ is the projection of $(\Thetabar_k, \Ybar_k)$ on $\Gamma_k$ and therefore $\sigma(\Thetahat_k,\Yhat_k) \subset \sigma(\Thetabar_k, \Ybar_k)$. Moreover, we have $\left\{(\Thetabar_{k},\Ybar_{k})=\gamma_k^p \right\} = \left\{(\Thetabar_{k},\Ybar_{k})=(\Thetahat_{k},\Yhat_{k})=\gamma_k^p \right\}$ so that  
 \begin{align*}
 |A_2|
   & = \big|\E[\Rbar'_k h(\Thetabar_{k}) - \Rbar'_k h(\Thetahat_{k}) \big| (\Thetahat_{k},\Yhat_{k})=\gamma_k^p]\big|\\
  &\leq (\barkpsup{h} + \fFtilde \overline{f}\barkplip{h}) \E[\dk{\Thetabar_{k} -\Thetahat_{k}} \big| (\Thetahat_{k},\Yhat_{k})=\gamma_k^p],
 \end{align*} 
 using lemma \ref{Rk'-lip}.\\
 $\bullet$  Finally for the last term, we have once again $A_3=0$.\\
 Taking the expectation yields the result.
\hfill$\Box$\\

\noindent\textbf{Proof Theorem \ref{lem:vbark-vhatk}}
One has
   \begin{eqnarray*}
     |\vhat'_N(\thetahat_N, \yhat) - \vbar_N'(\thetabar_N,\ybar)|    &\leq& |\min_{a\in\mathcal{M}} C'(\thetahat_N,\yhat,a) - \min_{a\in\mathcal{M}} C'(\thetabar_N,\ybar,a)|,
   \end{eqnarray*} 
 and we conclude using  Lemma \ref{lem:ga-BL}.
  
  Let $k\leq N-1$. Let $\thetabar_k$ on $\Omega_k$, $\ybar \in \Y$ and $(\thetahat_k,\yhat)=\gamma \in \Gamma_k$.
  \begin{align*}
    \lefteqn{|\vhat'_k(\thetahat_k, \yhat) - \vbar_k'(\thetabar_k,\ybar)|}\\
     & \leq    |\vbar'_k(\thetahat_k, \yhat) - \vbar_k'(\thetabar_k,\ybar)| +    |\vbar_k'(\thetahat_k,\yhat) - \vhat'_k(\thetahat_k, \yhat) |
  \end{align*}
  The first term is bounded by $\barklip{\vbar'_k} \dk{\thetahat_k-\thetabar_k}$ using Lemma \ref{stability'}. The second term is further decomposed into
  \begin{align*}
    \left| \vbar'_k(\thetahat_k, \yhat) - \vhat_k'(\thetahat_k,\yhat) \right| & \leq \left|\Rbar'_k\vbar'_{k+1}(\thetahat_k,\yhat)-\Rhat'_k\vhat'_{k+1}(\thetahat_k,\yhat)\right|  \\
    & \leq   \left|\Rhat'_k\vbar'_{k+1}(\thetahat_k,\yhat)-\Rhat'_k\vhat'_{k+1}(\thetahat_k,\yhat)\right|  \qquad (i)\\
    & \quad +   \left|\Rhat'_k\vbar'_{k+1}(\thetahat_k,\yhat)-\Rbar'_k\vbar'_{k+1}(\thetahat_k,\yhat)\right|  \qquad (ii)
  \end{align*}
  \begin{enumerate}
  \item[(i)]We have $ \left|\Rhat'_k\vhat'_{k+1}(\thetahat_k,\yhat)-\Rhat'_k\vbar'_{k+1}(\thetahat_k,\yhat)\right|   = \left|\Rhat'_k[\vhat'_{k+1} -\vbar'_{k+1}](\thetahat_k,\yhat)\right|$ and since $\Rhat'_k$ is a Markov kernel, \\ we have $ \left|\Rhat'_k\vhat'_{k+1}(\thetahat_k,\yhat)-\Rhat'_k\vbar'_{k+1}(\thetahat_k,\yhat)\right|\leq \barkpsup{\vhat'_{k+1}-\vbar'_{k+1}}$ where the supremum is taken over $\Gamma_{k+1}$.
   
  \item[(ii)] We use the second statement of Lemma \ref{lem:Rk-Rkbar} to get 
    \begin{align*}
      \lefteqn{ \left|\Rhat'_k\vbar'_{k+1}(\thetahat_k,\yhat)-\Rbar'_k\vbar'_{k+1}(\thetahat_k,\yhat)\right|}\\
      & \leq  \barkplip{\vbar'_{k+1}}\E\left[\dkp{\Thetahat_{k+1}-\Thetabar_{k+1}}\right]+ \left(\barkpsup{\vbar'_{k+1}}+\overline{f}\fFtilde \barkplip{\vbar'_{k+1}}\right) \E\left[\dk{\Thetahat_{k}-\Thetabar_{k}}\right]
    \end{align*}
  \end{enumerate}
  Hence combining the previous results leads to the expected result.
  \hfill$\Box$

\end{document}